\documentclass[a4paper,12pt]{article}

\usepackage{amsmath}

\usepackage{amssymb}
\usepackage{amsthm}
\usepackage[english]{babel}
\usepackage[latin1]{inputenc}
\usepackage{color}
\usepackage{graphicx}
\usepackage{xspace}
\usepackage{enumerate}
\usepackage{hyperref}
\hypersetup{hidelinks}
\usepackage[all]{xy}
\renewcommand{\phi}{\varphi}
\newtheorem{theor+}{Theorem}[section]
\newtheorem{propo+}[theor+]{Proposition}
\newtheorem{estim+}[theor+]{Estimate}
\newtheorem{imple+}[theor+]{Implementation remark}
\newtheorem{corol+}[theor+]{Corollary}
\newtheorem{notat+}[theor+]{Notation}
\newtheorem{lemma+}[theor+]{Lemma}
\newtheorem{obser+}[theor+]{Observation}
\theoremstyle{definition}
\newtheorem{defin+}[theor+]{Definition}
\newtheorem{examp+}[theor+]{Example}

\newtheorem{remar+}[theor+]{Remark}
\newtheorem{probl+}[theor+]{Problem}
\theoremstyle{remark}
\newtheorem{quest+}[theor+]{Question}
\newtheorem{pictu+}[theor+]{Picture}
\newtheorem{notes+}[theor+]{Notes}
\newtheorem{exerc+}[theor+]{Exercise}
\newtheorem{backg+}[theor+]{Background}
\newtheorem{fact+}[theor+]{Fact}
\newtheorem{abort+}[theor+]{Terminal state}

\newenvironment{theor}[1]{\begin{theor+}\label{#1}\slshape}{\end{theor+}}
\newenvironment{propo}[1]{\begin{propo+}\label{#1}\slshape}{\end{propo+}}

\newenvironment{lemma}[1]{\begin{lemma+}\label{#1}\slshape}{\end{lemma+}}

\newenvironment{defin}[1]{\begin{defin+}\label{#1}\upshape}{\end{defin+}}

\newenvironment{remar}[1]{\begin{remar+}\label{#1}\upshape}{\end{remar+}}
\newenvironment{probl}[1]{\begin{probl+}\label{#1}\upshape}{\end{probl+}}

\usepackage{pgfplots}
\pgfplotsset{compat=1.16}
\usetikzlibrary{arrows.meta}

\definecolor{expcol}{RGB}{38,110,200}
\definecolor{dircol}{RGB}{80,80,80}

\newcommand{\Z}{\mathbb Z}
\newcommand{\Aut}{\operatorname{Aut}}

\usepackage{stmaryrd}
\newcommand{\Iver}[1]{\llbracket #1\rrbracket}
\numberwithin{equation}{section}

\newcommand{\Ano}[1]{\cite[{\sf A#1}]{oeis}}

\newcommand{\lego}{\textsf{LEGO}\xspace}

\newcommand{\RR}{\mathbb{R}}

\newcommand{\ZZ}{\mathbb{Z}}
\newcommand{\QQ}{\mathbb{Q}}

\newcommand{\En}{E([n])}
\newcommand{\comp}{\operatorname{comp}}

\usepackage{pgfplots}
\pgfplotsset{compat=1.16}

\title{Efficient counting of  $2\times 4$ \lego\ brick buildings}
\author{S\o{}ren Eilers\thanks{Department of
    Mathematical Sciences, University
  of Copenhagen, Universitetsparken 5, DK-2100 Copenhagen \O, Denmark}}
\date{November 4, 2026}
\begin{document}
\maketitle

\begin{abstract}
We outline the algorithmic advances behind our recent computation of the number 
\[
8364925094730088313
\]
of ways to combine eleven $2\times 4$ \lego bricks.
\end{abstract}

\section{Background}
The number of ways to
combine six of the \lego\ Company's
iconic $2\times 4$  bricks is an
important figure for the company and
for enthusiasts alike. As detailed in \cite{legocprob}, it was obtained
by the author in 2004 after noticing that
a number previously used by the
company was too low.  

\begin{defin}{def:count}
We denote by $a(n)$ the number of ways to combine $n$ $2\times 4$ bricks placed parallel to the axes, in such a way that the \emph{contact graph} describing attached bricks is connected, identifying buildings up to translations in $\RR^3$ and rotations in the $XY$-plane.
\end{defin}

\begin{figure}
\begin{center}
\definecolor{blueT}{RGB}{40,104,200}
\definecolor{blueX}{RGB}{24,72,150}
\definecolor{blueY}{RGB}{32,90,180}
\definecolor{greenT}{RGB}{60,160,80}
\definecolor{greenX}{RGB}{36,116,54}
\definecolor{greenY}{RGB}{50,142,68}
\definecolor{redT}{RGB}{216,46,46}
\definecolor{redX}{RGB}{160,26,26}
\definecolor{redY}{RGB}{190,36,36}
\definecolor{whiteT}{RGB}{242,242,238}
\definecolor{whiteX}{RGB}{192,192,188}
\definecolor{whiteY}{RGB}{222,222,218}
\definecolor{yellowT}{RGB}{250,205,40}
\definecolor{yellowX}{RGB}{204,160,18}
\definecolor{yellowY}{RGB}{232,186,30}

\begin{tikzpicture}[scale=0.55]
\filldraw[fill=redX, draw=black!75, line width=0.5pt, line join=round] (4.0000,0.0000) -- (4.8400,0.6400) -- (4.8400,1.8400) -- (4.0000,1.2000) -- cycle;
\filldraw[fill=redY, draw=black!75, line width=0.5pt, line join=round] (0.0000,0.0000) -- (4.0000,0.0000) -- (4.0000,1.2000) -- (0.0000,1.2000) -- cycle;
\filldraw[fill=redT, draw=black!75, line width=0.5pt, line join=round] (0.0000,1.2000) -- (4.0000,1.2000) -- (4.8400,1.8400) -- (0.8400,1.8400) -- cycle;
\filldraw[fill=redY, draw=black!75, line width=0.3pt] (0.8046,1.6428) -- (0.8157,1.6212) -- (0.8482,1.6036) -- (0.8999,1.5911) -- (0.9673,1.5848) -- (1.0458,1.5849) -- (1.1300,1.5915) -- (1.2142,1.6041) -- (1.2927,1.6219) -- (1.3601,1.6437) -- (1.4118,1.6679) -- (1.4443,1.6930) -- (1.4554,1.7172) -- (1.4554,1.9172) -- (1.4443,1.9388) -- (1.4118,1.9564) -- (1.3601,1.9689) -- (1.2927,1.9752) -- (1.2142,1.9751) -- (1.1300,1.9685) -- (1.0458,1.9559) -- (0.9673,1.9381) -- (0.8999,1.9163) -- (0.8482,1.8921) -- (0.8157,1.8670) -- (0.8046,1.8428) -- cycle;
\filldraw[fill=redT, draw=black!75, line width=0.3pt] (1.4300,1.8800) -- (1.4524,1.9048) -- (1.4528,1.9280) -- (1.4312,1.9479) -- (1.3891,1.9631) -- (1.3294,1.9727) -- (1.2560,1.9760) -- (1.1741,1.9727) -- (1.0891,1.9631) -- (1.0070,1.9479) -- (0.9332,1.9280) -- (0.8728,1.9048) -- (0.8300,1.8800) -- (0.8076,1.8552) -- (0.8072,1.8320) -- (0.8288,1.8121) -- (0.8709,1.7969) -- (0.9306,1.7873) -- (1.0040,1.7840) -- (1.0859,1.7873) -- (1.1709,1.7969) -- (1.2530,1.8121) -- (1.3268,1.8320) -- (1.3872,1.8552) -- cycle;
\filldraw[fill=redY, draw=black!75, line width=0.3pt] (1.8046,1.6428) -- (1.8157,1.6212) -- (1.8482,1.6036) -- (1.8999,1.5911) -- (1.9673,1.5848) -- (2.0458,1.5849) -- (2.1300,1.5915) -- (2.2142,1.6041) -- (2.2927,1.6219) -- (2.3601,1.6437) -- (2.4118,1.6679) -- (2.4443,1.6930) -- (2.4554,1.7172) -- (2.4554,1.9172) -- (2.4443,1.9388) -- (2.4118,1.9564) -- (2.3601,1.9689) -- (2.2927,1.9752) -- (2.2142,1.9751) -- (2.1300,1.9685) -- (2.0458,1.9559) -- (1.9673,1.9381) -- (1.8999,1.9163) -- (1.8482,1.8921) -- (1.8157,1.8670) -- (1.8046,1.8428) -- cycle;
\filldraw[fill=redT, draw=black!75, line width=0.3pt] (2.4300,1.8800) -- (2.4524,1.9048) -- (2.4528,1.9280) -- (2.4312,1.9479) -- (2.3891,1.9631) -- (2.3294,1.9727) -- (2.2560,1.9760) -- (2.1741,1.9727) -- (2.0891,1.9631) -- (2.0070,1.9479) -- (1.9332,1.9280) -- (1.8728,1.9048) -- (1.8300,1.8800) -- (1.8076,1.8552) -- (1.8072,1.8320) -- (1.8288,1.8121) -- (1.8709,1.7969) -- (1.9306,1.7873) -- (2.0040,1.7840) -- (2.0859,1.7873) -- (2.1709,1.7969) -- (2.2530,1.8121) -- (2.3268,1.8320) -- (2.3872,1.8552) -- cycle;
\filldraw[fill=redY, draw=black!75, line width=0.3pt] (2.8046,1.6428) -- (2.8157,1.6212) -- (2.8482,1.6036) -- (2.8999,1.5911) -- (2.9673,1.5848) -- (3.0458,1.5849) -- (3.1300,1.5915) -- (3.2142,1.6041) -- (3.2927,1.6219) -- (3.3601,1.6437) -- (3.4118,1.6679) -- (3.4443,1.6930) -- (3.4554,1.7172) -- (3.4554,1.9172) -- (3.4443,1.9388) -- (3.4118,1.9564) -- (3.3601,1.9689) -- (3.2927,1.9752) -- (3.2142,1.9751) -- (3.1300,1.9685) -- (3.0458,1.9559) -- (2.9673,1.9381) -- (2.8999,1.9163) -- (2.8482,1.8921) -- (2.8157,1.8670) -- (2.8046,1.8428) -- cycle;
\filldraw[fill=redT, draw=black!75, line width=0.3pt] (3.4300,1.8800) -- (3.4524,1.9048) -- (3.4528,1.9280) -- (3.4312,1.9479) -- (3.3891,1.9631) -- (3.3294,1.9727) -- (3.2560,1.9760) -- (3.1741,1.9727) -- (3.0891,1.9631) -- (3.0070,1.9479) -- (2.9332,1.9280) -- (2.8728,1.9048) -- (2.8300,1.8800) -- (2.8076,1.8552) -- (2.8072,1.8320) -- (2.8288,1.8121) -- (2.8709,1.7969) -- (2.9306,1.7873) -- (3.0040,1.7840) -- (3.0859,1.7873) -- (3.1709,1.7969) -- (3.2530,1.8121) -- (3.3268,1.8320) -- (3.3872,1.8552) -- cycle;
\filldraw[fill=redY, draw=black!75, line width=0.3pt] (3.8046,1.6428) -- (3.8157,1.6212) -- (3.8482,1.6036) -- (3.8999,1.5911) -- (3.9673,1.5848) -- (4.0458,1.5849) -- (4.1300,1.5915) -- (4.2142,1.6041) -- (4.2927,1.6219) -- (4.3601,1.6437) -- (4.4118,1.6679) -- (4.4443,1.6930) -- (4.4554,1.7172) -- (4.4554,1.9172) -- (4.4443,1.9388) -- (4.4118,1.9564) -- (4.3601,1.9689) -- (4.2927,1.9752) -- (4.2142,1.9751) -- (4.1300,1.9685) -- (4.0458,1.9559) -- (3.9673,1.9381) -- (3.8999,1.9163) -- (3.8482,1.8921) -- (3.8157,1.8670) -- (3.8046,1.8428) -- cycle;
\filldraw[fill=redT, draw=black!75, line width=0.3pt] (4.4300,1.8800) -- (4.4524,1.9048) -- (4.4528,1.9280) -- (4.4312,1.9479) -- (4.3891,1.9631) -- (4.3294,1.9727) -- (4.2560,1.9760) -- (4.1741,1.9727) -- (4.0891,1.9631) -- (4.0070,1.9479) -- (3.9332,1.9280) -- (3.8728,1.9048) -- (3.8300,1.8800) -- (3.8076,1.8552) -- (3.8072,1.8320) -- (3.8288,1.8121) -- (3.8709,1.7969) -- (3.9306,1.7873) -- (4.0040,1.7840) -- (4.0859,1.7873) -- (4.1709,1.7969) -- (4.2530,1.8121) -- (4.3268,1.8320) -- (4.3872,1.8552) -- cycle;
\filldraw[fill=redY, draw=black!75, line width=0.3pt] (0.3846,1.3228) -- (0.3957,1.3012) -- (0.4282,1.2836) -- (0.4799,1.2711) -- (0.5473,1.2648) -- (0.6258,1.2649) -- (0.7100,1.2715) -- (0.7942,1.2841) -- (0.8727,1.3019) -- (0.9401,1.3237) -- (0.9918,1.3479) -- (1.0243,1.3730) -- (1.0354,1.3972) -- (1.0354,1.5972) -- (1.0243,1.6188) -- (0.9918,1.6364) -- (0.9401,1.6489) -- (0.8727,1.6552) -- (0.7942,1.6551) -- (0.7100,1.6485) -- (0.6258,1.6359) -- (0.5473,1.6181) -- (0.4799,1.5963) -- (0.4282,1.5721) -- (0.3957,1.5470) -- (0.3846,1.5228) -- cycle;
\filldraw[fill=redT, draw=black!75, line width=0.3pt] (1.0100,1.5600) -- (1.0324,1.5848) -- (1.0328,1.6080) -- (1.0112,1.6279) -- (0.9691,1.6431) -- (0.9094,1.6527) -- (0.8360,1.6560) -- (0.7541,1.6527) -- (0.6691,1.6431) -- (0.5870,1.6279) -- (0.5132,1.6080) -- (0.4528,1.5848) -- (0.4100,1.5600) -- (0.3876,1.5352) -- (0.3872,1.5120) -- (0.4088,1.4921) -- (0.4509,1.4769) -- (0.5106,1.4673) -- (0.5840,1.4640) -- (0.6659,1.4673) -- (0.7509,1.4769) -- (0.8330,1.4921) -- (0.9068,1.5120) -- (0.9672,1.5352) -- cycle;
\filldraw[fill=redY, draw=black!75, line width=0.3pt] (1.3846,1.3228) -- (1.3957,1.3012) -- (1.4282,1.2836) -- (1.4799,1.2711) -- (1.5473,1.2648) -- (1.6258,1.2649) -- (1.7100,1.2715) -- (1.7942,1.2841) -- (1.8727,1.3019) -- (1.9401,1.3237) -- (1.9918,1.3479) -- (2.0243,1.3730) -- (2.0354,1.3972) -- (2.0354,1.5972) -- (2.0243,1.6188) -- (1.9918,1.6364) -- (1.9401,1.6489) -- (1.8727,1.6552) -- (1.7942,1.6551) -- (1.7100,1.6485) -- (1.6258,1.6359) -- (1.5473,1.6181) -- (1.4799,1.5963) -- (1.4282,1.5721) -- (1.3957,1.5470) -- (1.3846,1.5228) -- cycle;
\filldraw[fill=redT, draw=black!75, line width=0.3pt] (2.0100,1.5600) -- (2.0324,1.5848) -- (2.0328,1.6080) -- (2.0112,1.6279) -- (1.9691,1.6431) -- (1.9094,1.6527) -- (1.8360,1.6560) -- (1.7541,1.6527) -- (1.6691,1.6431) -- (1.5870,1.6279) -- (1.5132,1.6080) -- (1.4528,1.5848) -- (1.4100,1.5600) -- (1.3876,1.5352) -- (1.3872,1.5120) -- (1.4088,1.4921) -- (1.4509,1.4769) -- (1.5106,1.4673) -- (1.5840,1.4640) -- (1.6659,1.4673) -- (1.7509,1.4769) -- (1.8330,1.4921) -- (1.9068,1.5120) -- (1.9672,1.5352) -- cycle;
\filldraw[fill=redY, draw=black!75, line width=0.3pt] (2.3846,1.3228) -- (2.3957,1.3012) -- (2.4282,1.2836) -- (2.4799,1.2711) -- (2.5473,1.2648) -- (2.6258,1.2649) -- (2.7100,1.2715) -- (2.7942,1.2841) -- (2.8727,1.3019) -- (2.9401,1.3237) -- (2.9918,1.3479) -- (3.0243,1.3730) -- (3.0354,1.3972) -- (3.0354,1.5972) -- (3.0243,1.6188) -- (2.9918,1.6364) -- (2.9401,1.6489) -- (2.8727,1.6552) -- (2.7942,1.6551) -- (2.7100,1.6485) -- (2.6258,1.6359) -- (2.5473,1.6181) -- (2.4799,1.5963) -- (2.4282,1.5721) -- (2.3957,1.5470) -- (2.3846,1.5228) -- cycle;
\filldraw[fill=redT, draw=black!75, line width=0.3pt] (3.0100,1.5600) -- (3.0324,1.5848) -- (3.0328,1.6080) -- (3.0112,1.6279) -- (2.9691,1.6431) -- (2.9094,1.6527) -- (2.8360,1.6560) -- (2.7541,1.6527) -- (2.6691,1.6431) -- (2.5870,1.6279) -- (2.5132,1.6080) -- (2.4528,1.5848) -- (2.4100,1.5600) -- (2.3876,1.5352) -- (2.3872,1.5120) -- (2.4088,1.4921) -- (2.4509,1.4769) -- (2.5106,1.4673) -- (2.5840,1.4640) -- (2.6659,1.4673) -- (2.7509,1.4769) -- (2.8330,1.4921) -- (2.9068,1.5120) -- (2.9672,1.5352) -- cycle;
\filldraw[fill=redY, draw=black!75, line width=0.3pt] (3.3846,1.3228) -- (3.3957,1.3012) -- (3.4282,1.2836) -- (3.4799,1.2711) -- (3.5473,1.2648) -- (3.6258,1.2649) -- (3.7100,1.2715) -- (3.7942,1.2841) -- (3.8727,1.3019) -- (3.9401,1.3237) -- (3.9918,1.3479) -- (4.0243,1.3730) -- (4.0354,1.3972) -- (4.0354,1.5972) -- (4.0243,1.6188) -- (3.9918,1.6364) -- (3.9401,1.6489) -- (3.8727,1.6552) -- (3.7942,1.6551) -- (3.7100,1.6485) -- (3.6258,1.6359) -- (3.5473,1.6181) -- (3.4799,1.5963) -- (3.4282,1.5721) -- (3.3957,1.5470) -- (3.3846,1.5228) -- cycle;
\filldraw[fill=redT, draw=black!75, line width=0.3pt] (4.0100,1.5600) -- (4.0324,1.5848) -- (4.0328,1.6080) -- (4.0112,1.6279) -- (3.9691,1.6431) -- (3.9094,1.6527) -- (3.8360,1.6560) -- (3.7541,1.6527) -- (3.6691,1.6431) -- (3.5870,1.6279) -- (3.5132,1.6080) -- (3.4528,1.5848) -- (3.4100,1.5600) -- (3.3876,1.5352) -- (3.3872,1.5120) -- (3.4088,1.4921) -- (3.4509,1.4769) -- (3.5106,1.4673) -- (3.5840,1.4640) -- (3.6659,1.4673) -- (3.7509,1.4769) -- (3.8330,1.4921) -- (3.9068,1.5120) -- (3.9672,1.5352) -- cycle;
\filldraw[fill=blueX, draw=black!75, line width=0.5pt, line join=round] (10.0000,0.0000) -- (10.8400,0.6400) -- (10.8400,1.8400) -- (10.0000,1.2000) -- cycle;
\filldraw[fill=blueY, draw=black!75, line width=0.5pt, line join=round] (6.0000,0.0000) -- (10.0000,0.0000) -- (10.0000,1.2000) -- (6.0000,1.2000) -- cycle;
\filldraw[fill=blueT, draw=black!75, line width=0.5pt, line join=round] (6.0000,1.2000) -- (10.0000,1.2000) -- (10.8400,1.8400) -- (6.8400,1.8400) -- cycle;
\filldraw[fill=blueY, draw=black!75, line width=0.3pt] (6.8046,1.6428) -- (6.8157,1.6212) -- (6.8482,1.6036) -- (6.8999,1.5911) -- (6.9673,1.5848) -- (7.0458,1.5849) -- (7.1300,1.5915) -- (7.2142,1.6041) -- (7.2927,1.6219) -- (7.3601,1.6437) -- (7.4118,1.6679) -- (7.4443,1.6930) -- (7.4554,1.7172) -- (7.4554,1.9172) -- (7.4443,1.9388) -- (7.4118,1.9564) -- (7.3601,1.9689) -- (7.2927,1.9752) -- (7.2142,1.9751) -- (7.1300,1.9685) -- (7.0458,1.9559) -- (6.9673,1.9381) -- (6.8999,1.9163) -- (6.8482,1.8921) -- (6.8157,1.8670) -- (6.8046,1.8428) -- cycle;
\filldraw[fill=blueT, draw=black!75, line width=0.3pt] (7.4300,1.8800) -- (7.4524,1.9048) -- (7.4528,1.9280) -- (7.4312,1.9479) -- (7.3891,1.9631) -- (7.3294,1.9727) -- (7.2560,1.9760) -- (7.1741,1.9727) -- (7.0891,1.9631) -- (7.0070,1.9479) -- (6.9332,1.9280) -- (6.8728,1.9048) -- (6.8300,1.8800) -- (6.8076,1.8552) -- (6.8072,1.8320) -- (6.8288,1.8121) -- (6.8709,1.7969) -- (6.9306,1.7873) -- (7.0040,1.7840) -- (7.0859,1.7873) -- (7.1709,1.7969) -- (7.2530,1.8121) -- (7.3268,1.8320) -- (7.3872,1.8552) -- cycle;
\filldraw[fill=blueY, draw=black!75, line width=0.3pt] (7.8046,1.6428) -- (7.8157,1.6212) -- (7.8482,1.6036) -- (7.8999,1.5911) -- (7.9673,1.5848) -- (8.0458,1.5849) -- (8.1300,1.5915) -- (8.2142,1.6041) -- (8.2927,1.6219) -- (8.3601,1.6437) -- (8.4118,1.6679) -- (8.4443,1.6930) -- (8.4554,1.7172) -- (8.4554,1.9172) -- (8.4443,1.9388) -- (8.4118,1.9564) -- (8.3601,1.9689) -- (8.2927,1.9752) -- (8.2142,1.9751) -- (8.1300,1.9685) -- (8.0458,1.9559) -- (7.9673,1.9381) -- (7.8999,1.9163) -- (7.8482,1.8921) -- (7.8157,1.8670) -- (7.8046,1.8428) -- cycle;
\filldraw[fill=blueT, draw=black!75, line width=0.3pt] (8.4300,1.8800) -- (8.4524,1.9048) -- (8.4528,1.9280) -- (8.4312,1.9479) -- (8.3891,1.9631) -- (8.3294,1.9727) -- (8.2560,1.9760) -- (8.1741,1.9727) -- (8.0891,1.9631) -- (8.0070,1.9479) -- (7.9332,1.9280) -- (7.8728,1.9048) -- (7.8300,1.8800) -- (7.8076,1.8552) -- (7.8072,1.8320) -- (7.8288,1.8121) -- (7.8709,1.7969) -- (7.9306,1.7873) -- (8.0040,1.7840) -- (8.0859,1.7873) -- (8.1709,1.7969) -- (8.2530,1.8121) -- (8.3268,1.8320) -- (8.3872,1.8552) -- cycle;
\filldraw[fill=blueY, draw=black!75, line width=0.3pt] (8.8046,1.6428) -- (8.8157,1.6212) -- (8.8482,1.6036) -- (8.8999,1.5911) -- (8.9673,1.5848) -- (9.0458,1.5849) -- (9.1300,1.5915) -- (9.2142,1.6041) -- (9.2927,1.6219) -- (9.3601,1.6437) -- (9.4118,1.6679) -- (9.4443,1.6930) -- (9.4554,1.7172) -- (9.4554,1.9172) -- (9.4443,1.9388) -- (9.4118,1.9564) -- (9.3601,1.9689) -- (9.2927,1.9752) -- (9.2142,1.9751) -- (9.1300,1.9685) -- (9.0458,1.9559) -- (8.9673,1.9381) -- (8.8999,1.9163) -- (8.8482,1.8921) -- (8.8157,1.8670) -- (8.8046,1.8428) -- cycle;
\filldraw[fill=blueT, draw=black!75, line width=0.3pt] (9.4300,1.8800) -- (9.4524,1.9048) -- (9.4528,1.9280) -- (9.4312,1.9479) -- (9.3891,1.9631) -- (9.3294,1.9727) -- (9.2560,1.9760) -- (9.1741,1.9727) -- (9.0891,1.9631) -- (9.0070,1.9479) -- (8.9332,1.9280) -- (8.8728,1.9048) -- (8.8300,1.8800) -- (8.8076,1.8552) -- (8.8072,1.8320) -- (8.8288,1.8121) -- (8.8709,1.7969) -- (8.9306,1.7873) -- (9.0040,1.7840) -- (9.0859,1.7873) -- (9.1709,1.7969) -- (9.2530,1.8121) -- (9.3268,1.8320) -- (9.3872,1.8552) -- cycle;
\filldraw[fill=blueY, draw=black!75, line width=0.3pt] (9.8046,1.6428) -- (9.8157,1.6212) -- (9.8482,1.6036) -- (9.8999,1.5911) -- (9.9673,1.5848) -- (10.0458,1.5849) -- (10.1300,1.5915) -- (10.2142,1.6041) -- (10.2927,1.6219) -- (10.3601,1.6437) -- (10.4118,1.6679) -- (10.4443,1.6930) -- (10.4554,1.7172) -- (10.4554,1.9172) -- (10.4443,1.9388) -- (10.4118,1.9564) -- (10.3601,1.9689) -- (10.2927,1.9752) -- (10.2142,1.9751) -- (10.1300,1.9685) -- (10.0458,1.9559) -- (9.9673,1.9381) -- (9.8999,1.9163) -- (9.8482,1.8921) -- (9.8157,1.8670) -- (9.8046,1.8428) -- cycle;
\filldraw[fill=blueT, draw=black!75, line width=0.3pt] (10.4300,1.8800) -- (10.4524,1.9048) -- (10.4528,1.9280) -- (10.4312,1.9479) -- (10.3891,1.9631) -- (10.3294,1.9727) -- (10.2560,1.9760) -- (10.1741,1.9727) -- (10.0891,1.9631) -- (10.0070,1.9479) -- (9.9332,1.9280) -- (9.8728,1.9048) -- (9.8300,1.8800) -- (9.8076,1.8552) -- (9.8072,1.8320) -- (9.8288,1.8121) -- (9.8709,1.7969) -- (9.9306,1.7873) -- (10.0040,1.7840) -- (10.0859,1.7873) -- (10.1709,1.7969) -- (10.2530,1.8121) -- (10.3268,1.8320) -- (10.3872,1.8552) -- cycle;
\filldraw[fill=blueY, draw=black!75, line width=0.3pt] (6.3846,1.3228) -- (6.3957,1.3012) -- (6.4282,1.2836) -- (6.4799,1.2711) -- (6.5473,1.2648) -- (6.6258,1.2649) -- (6.7100,1.2715) -- (6.7942,1.2841) -- (6.8727,1.3019) -- (6.9401,1.3237) -- (6.9918,1.3479) -- (7.0243,1.3730) -- (7.0354,1.3972) -- (7.0354,1.5972) -- (7.0243,1.6188) -- (6.9918,1.6364) -- (6.9401,1.6489) -- (6.8727,1.6552) -- (6.7942,1.6551) -- (6.7100,1.6485) -- (6.6258,1.6359) -- (6.5473,1.6181) -- (6.4799,1.5963) -- (6.4282,1.5721) -- (6.3957,1.5470) -- (6.3846,1.5228) -- cycle;
\filldraw[fill=blueT, draw=black!75, line width=0.3pt] (7.0100,1.5600) -- (7.0324,1.5848) -- (7.0328,1.6080) -- (7.0112,1.6279) -- (6.9691,1.6431) -- (6.9094,1.6527) -- (6.8360,1.6560) -- (6.7541,1.6527) -- (6.6691,1.6431) -- (6.5870,1.6279) -- (6.5132,1.6080) -- (6.4528,1.5848) -- (6.4100,1.5600) -- (6.3876,1.5352) -- (6.3872,1.5120) -- (6.4088,1.4921) -- (6.4509,1.4769) -- (6.5106,1.4673) -- (6.5840,1.4640) -- (6.6659,1.4673) -- (6.7509,1.4769) -- (6.8330,1.4921) -- (6.9068,1.5120) -- (6.9672,1.5352) -- cycle;
\filldraw[fill=blueY, draw=black!75, line width=0.3pt] (7.3846,1.3228) -- (7.3957,1.3012) -- (7.4282,1.2836) -- (7.4799,1.2711) -- (7.5473,1.2648) -- (7.6258,1.2649) -- (7.7100,1.2715) -- (7.7942,1.2841) -- (7.8727,1.3019) -- (7.9401,1.3237) -- (7.9918,1.3479) -- (8.0243,1.3730) -- (8.0354,1.3972) -- (8.0354,1.5972) -- (8.0243,1.6188) -- (7.9918,1.6364) -- (7.9401,1.6489) -- (7.8727,1.6552) -- (7.7942,1.6551) -- (7.7100,1.6485) -- (7.6258,1.6359) -- (7.5473,1.6181) -- (7.4799,1.5963) -- (7.4282,1.5721) -- (7.3957,1.5470) -- (7.3846,1.5228) -- cycle;
\filldraw[fill=blueT, draw=black!75, line width=0.3pt] (8.0100,1.5600) -- (8.0324,1.5848) -- (8.0328,1.6080) -- (8.0112,1.6279) -- (7.9691,1.6431) -- (7.9094,1.6527) -- (7.8360,1.6560) -- (7.7541,1.6527) -- (7.6691,1.6431) -- (7.5870,1.6279) -- (7.5132,1.6080) -- (7.4528,1.5848) -- (7.4100,1.5600) -- (7.3876,1.5352) -- (7.3872,1.5120) -- (7.4088,1.4921) -- (7.4509,1.4769) -- (7.5106,1.4673) -- (7.5840,1.4640) -- (7.6659,1.4673) -- (7.7509,1.4769) -- (7.8330,1.4921) -- (7.9068,1.5120) -- (7.9672,1.5352) -- cycle;
\filldraw[fill=blueY, draw=black!75, line width=0.3pt] (8.3846,1.3228) -- (8.3957,1.3012) -- (8.4282,1.2836) -- (8.4799,1.2711) -- (8.5473,1.2648) -- (8.6258,1.2649) -- (8.7100,1.2715) -- (8.7942,1.2841) -- (8.8727,1.3019) -- (8.9401,1.3237) -- (8.9918,1.3479) -- (9.0243,1.3730) -- (9.0354,1.3972) -- (9.0354,1.5972) -- (9.0243,1.6188) -- (8.9918,1.6364) -- (8.9401,1.6489) -- (8.8727,1.6552) -- (8.7942,1.6551) -- (8.7100,1.6485) -- (8.6258,1.6359) -- (8.5473,1.6181) -- (8.4799,1.5963) -- (8.4282,1.5721) -- (8.3957,1.5470) -- (8.3846,1.5228) -- cycle;
\filldraw[fill=blueT, draw=black!75, line width=0.3pt] (9.0100,1.5600) -- (9.0324,1.5848) -- (9.0328,1.6080) -- (9.0112,1.6279) -- (8.9691,1.6431) -- (8.9094,1.6527) -- (8.8360,1.6560) -- (8.7541,1.6527) -- (8.6691,1.6431) -- (8.5870,1.6279) -- (8.5132,1.6080) -- (8.4528,1.5848) -- (8.4100,1.5600) -- (8.3876,1.5352) -- (8.3872,1.5120) -- (8.4088,1.4921) -- (8.4509,1.4769) -- (8.5106,1.4673) -- (8.5840,1.4640) -- (8.6659,1.4673) -- (8.7509,1.4769) -- (8.8330,1.4921) -- (8.9068,1.5120) -- (8.9672,1.5352) -- cycle;
\filldraw[fill=blueY, draw=black!75, line width=0.3pt] (9.3846,1.3228) -- (9.3957,1.3012) -- (9.4282,1.2836) -- (9.4799,1.2711) -- (9.5473,1.2648) -- (9.6258,1.2649) -- (9.7100,1.2715) -- (9.7942,1.2841) -- (9.8727,1.3019) -- (9.9401,1.3237) -- (9.9918,1.3479) -- (10.0243,1.3730) -- (10.0354,1.3972) -- (10.0354,1.5972) -- (10.0243,1.6188) -- (9.9918,1.6364) -- (9.9401,1.6489) -- (9.8727,1.6552) -- (9.7942,1.6551) -- (9.7100,1.6485) -- (9.6258,1.6359) -- (9.5473,1.6181) -- (9.4799,1.5963) -- (9.4282,1.5721) -- (9.3957,1.5470) -- (9.3846,1.5228) -- cycle;
\filldraw[fill=blueT, draw=black!75, line width=0.3pt] (10.0100,1.5600) -- (10.0324,1.5848) -- (10.0328,1.6080) -- (10.0112,1.6279) -- (9.9691,1.6431) -- (9.9094,1.6527) -- (9.8360,1.6560) -- (9.7541,1.6527) -- (9.6691,1.6431) -- (9.5870,1.6279) -- (9.5132,1.6080) -- (9.4528,1.5848) -- (9.4100,1.5600) -- (9.3876,1.5352) -- (9.3872,1.5120) -- (9.4088,1.4921) -- (9.4509,1.4769) -- (9.5106,1.4673) -- (9.5840,1.4640) -- (9.6659,1.4673) -- (9.7509,1.4769) -- (9.8330,1.4921) -- (9.9068,1.5120) -- (9.9672,1.5352) -- cycle;
\filldraw[fill=yellowX, draw=black!75, line width=0.5pt, line join=round] (7.0000,1.2000) -- (7.8400,1.8400) -- (7.8400,3.0400) -- (7.0000,2.4000) -- cycle;
\filldraw[fill=yellowY, draw=black!75, line width=0.5pt, line join=round] (3.0000,1.2000) -- (7.0000,1.2000) -- (7.0000,2.4000) -- (3.0000,2.4000) -- cycle;
\filldraw[fill=yellowT, draw=black!75, line width=0.5pt, line join=round] (3.0000,2.4000) -- (7.0000,2.4000) -- (7.8400,3.0400) -- (3.8400,3.0400) -- cycle;
\filldraw[fill=yellowY, draw=black!75, line width=0.3pt] (3.8046,2.8428) -- (3.8157,2.8212) -- (3.8482,2.8036) -- (3.8999,2.7911) -- (3.9673,2.7848) -- (4.0458,2.7849) -- (4.1300,2.7915) -- (4.2142,2.8041) -- (4.2927,2.8219) -- (4.3601,2.8437) -- (4.4118,2.8679) -- (4.4443,2.8930) -- (4.4554,2.9172) -- (4.4554,3.1172) -- (4.4443,3.1388) -- (4.4118,3.1564) -- (4.3601,3.1689) -- (4.2927,3.1752) -- (4.2142,3.1751) -- (4.1300,3.1685) -- (4.0458,3.1559) -- (3.9673,3.1381) -- (3.8999,3.1163) -- (3.8482,3.0921) -- (3.8157,3.0670) -- (3.8046,3.0428) -- cycle;
\filldraw[fill=yellowT, draw=black!75, line width=0.3pt] (4.4300,3.0800) -- (4.4524,3.1048) -- (4.4528,3.1280) -- (4.4312,3.1479) -- (4.3891,3.1631) -- (4.3294,3.1727) -- (4.2560,3.1760) -- (4.1741,3.1727) -- (4.0891,3.1631) -- (4.0070,3.1479) -- (3.9332,3.1280) -- (3.8728,3.1048) -- (3.8300,3.0800) -- (3.8076,3.0552) -- (3.8072,3.0320) -- (3.8288,3.0121) -- (3.8709,2.9969) -- (3.9306,2.9873) -- (4.0040,2.9840) -- (4.0859,2.9873) -- (4.1709,2.9969) -- (4.2530,3.0121) -- (4.3268,3.0320) -- (4.3872,3.0552) -- cycle;
\filldraw[fill=yellowY, draw=black!75, line width=0.3pt] (4.8046,2.8428) -- (4.8157,2.8212) -- (4.8482,2.8036) -- (4.8999,2.7911) -- (4.9673,2.7848) -- (5.0458,2.7849) -- (5.1300,2.7915) -- (5.2142,2.8041) -- (5.2927,2.8219) -- (5.3601,2.8437) -- (5.4118,2.8679) -- (5.4443,2.8930) -- (5.4554,2.9172) -- (5.4554,3.1172) -- (5.4443,3.1388) -- (5.4118,3.1564) -- (5.3601,3.1689) -- (5.2927,3.1752) -- (5.2142,3.1751) -- (5.1300,3.1685) -- (5.0458,3.1559) -- (4.9673,3.1381) -- (4.8999,3.1163) -- (4.8482,3.0921) -- (4.8157,3.0670) -- (4.8046,3.0428) -- cycle;
\filldraw[fill=yellowT, draw=black!75, line width=0.3pt] (5.4300,3.0800) -- (5.4524,3.1048) -- (5.4528,3.1280) -- (5.4312,3.1479) -- (5.3891,3.1631) -- (5.3294,3.1727) -- (5.2560,3.1760) -- (5.1741,3.1727) -- (5.0891,3.1631) -- (5.0070,3.1479) -- (4.9332,3.1280) -- (4.8728,3.1048) -- (4.8300,3.0800) -- (4.8076,3.0552) -- (4.8072,3.0320) -- (4.8288,3.0121) -- (4.8709,2.9969) -- (4.9306,2.9873) -- (5.0040,2.9840) -- (5.0859,2.9873) -- (5.1709,2.9969) -- (5.2530,3.0121) -- (5.3268,3.0320) -- (5.3872,3.0552) -- cycle;
\filldraw[fill=yellowY, draw=black!75, line width=0.3pt] (5.8046,2.8428) -- (5.8157,2.8212) -- (5.8482,2.8036) -- (5.8999,2.7911) -- (5.9673,2.7848) -- (6.0458,2.7849) -- (6.1300,2.7915) -- (6.2142,2.8041) -- (6.2927,2.8219) -- (6.3601,2.8437) -- (6.4118,2.8679) -- (6.4443,2.8930) -- (6.4554,2.9172) -- (6.4554,3.1172) -- (6.4443,3.1388) -- (6.4118,3.1564) -- (6.3601,3.1689) -- (6.2927,3.1752) -- (6.2142,3.1751) -- (6.1300,3.1685) -- (6.0458,3.1559) -- (5.9673,3.1381) -- (5.8999,3.1163) -- (5.8482,3.0921) -- (5.8157,3.0670) -- (5.8046,3.0428) -- cycle;
\filldraw[fill=yellowT, draw=black!75, line width=0.3pt] (6.4300,3.0800) -- (6.4524,3.1048) -- (6.4528,3.1280) -- (6.4312,3.1479) -- (6.3891,3.1631) -- (6.3294,3.1727) -- (6.2560,3.1760) -- (6.1741,3.1727) -- (6.0891,3.1631) -- (6.0070,3.1479) -- (5.9332,3.1280) -- (5.8728,3.1048) -- (5.8300,3.0800) -- (5.8076,3.0552) -- (5.8072,3.0320) -- (5.8288,3.0121) -- (5.8709,2.9969) -- (5.9306,2.9873) -- (6.0040,2.9840) -- (6.0859,2.9873) -- (6.1709,2.9969) -- (6.2530,3.0121) -- (6.3268,3.0320) -- (6.3872,3.0552) -- cycle;
\filldraw[fill=yellowY, draw=black!75, line width=0.3pt] (6.8046,2.8428) -- (6.8157,2.8212) -- (6.8482,2.8036) -- (6.8999,2.7911) -- (6.9673,2.7848) -- (7.0458,2.7849) -- (7.1300,2.7915) -- (7.2142,2.8041) -- (7.2927,2.8219) -- (7.3601,2.8437) -- (7.4118,2.8679) -- (7.4443,2.8930) -- (7.4554,2.9172) -- (7.4554,3.1172) -- (7.4443,3.1388) -- (7.4118,3.1564) -- (7.3601,3.1689) -- (7.2927,3.1752) -- (7.2142,3.1751) -- (7.1300,3.1685) -- (7.0458,3.1559) -- (6.9673,3.1381) -- (6.8999,3.1163) -- (6.8482,3.0921) -- (6.8157,3.0670) -- (6.8046,3.0428) -- cycle;
\filldraw[fill=yellowT, draw=black!75, line width=0.3pt] (7.4300,3.0800) -- (7.4524,3.1048) -- (7.4528,3.1280) -- (7.4312,3.1479) -- (7.3891,3.1631) -- (7.3294,3.1727) -- (7.2560,3.1760) -- (7.1741,3.1727) -- (7.0891,3.1631) -- (7.0070,3.1479) -- (6.9332,3.1280) -- (6.8728,3.1048) -- (6.8300,3.0800) -- (6.8076,3.0552) -- (6.8072,3.0320) -- (6.8288,3.0121) -- (6.8709,2.9969) -- (6.9306,2.9873) -- (7.0040,2.9840) -- (7.0859,2.9873) -- (7.1709,2.9969) -- (7.2530,3.0121) -- (7.3268,3.0320) -- (7.3872,3.0552) -- cycle;
\filldraw[fill=yellowY, draw=black!75, line width=0.3pt] (3.3846,2.5228) -- (3.3957,2.5012) -- (3.4282,2.4836) -- (3.4799,2.4711) -- (3.5473,2.4648) -- (3.6258,2.4649) -- (3.7100,2.4715) -- (3.7942,2.4841) -- (3.8727,2.5019) -- (3.9401,2.5237) -- (3.9918,2.5479) -- (4.0243,2.5730) -- (4.0354,2.5972) -- (4.0354,2.7972) -- (4.0243,2.8188) -- (3.9918,2.8364) -- (3.9401,2.8489) -- (3.8727,2.8552) -- (3.7942,2.8551) -- (3.7100,2.8485) -- (3.6258,2.8359) -- (3.5473,2.8181) -- (3.4799,2.7963) -- (3.4282,2.7721) -- (3.3957,2.7470) -- (3.3846,2.7228) -- cycle;
\filldraw[fill=yellowT, draw=black!75, line width=0.3pt] (4.0100,2.7600) -- (4.0324,2.7848) -- (4.0328,2.8080) -- (4.0112,2.8279) -- (3.9691,2.8431) -- (3.9094,2.8527) -- (3.8360,2.8560) -- (3.7541,2.8527) -- (3.6691,2.8431) -- (3.5870,2.8279) -- (3.5132,2.8080) -- (3.4528,2.7848) -- (3.4100,2.7600) -- (3.3876,2.7352) -- (3.3872,2.7120) -- (3.4088,2.6921) -- (3.4509,2.6769) -- (3.5106,2.6673) -- (3.5840,2.6640) -- (3.6659,2.6673) -- (3.7509,2.6769) -- (3.8330,2.6921) -- (3.9068,2.7120) -- (3.9672,2.7352) -- cycle;
\filldraw[fill=yellowY, draw=black!75, line width=0.3pt] (4.3846,2.5228) -- (4.3957,2.5012) -- (4.4282,2.4836) -- (4.4799,2.4711) -- (4.5473,2.4648) -- (4.6258,2.4649) -- (4.7100,2.4715) -- (4.7942,2.4841) -- (4.8727,2.5019) -- (4.9401,2.5237) -- (4.9918,2.5479) -- (5.0243,2.5730) -- (5.0354,2.5972) -- (5.0354,2.7972) -- (5.0243,2.8188) -- (4.9918,2.8364) -- (4.9401,2.8489) -- (4.8727,2.8552) -- (4.7942,2.8551) -- (4.7100,2.8485) -- (4.6258,2.8359) -- (4.5473,2.8181) -- (4.4799,2.7963) -- (4.4282,2.7721) -- (4.3957,2.7470) -- (4.3846,2.7228) -- cycle;
\filldraw[fill=yellowT, draw=black!75, line width=0.3pt] (5.0100,2.7600) -- (5.0324,2.7848) -- (5.0328,2.8080) -- (5.0112,2.8279) -- (4.9691,2.8431) -- (4.9094,2.8527) -- (4.8360,2.8560) -- (4.7541,2.8527) -- (4.6691,2.8431) -- (4.5870,2.8279) -- (4.5132,2.8080) -- (4.4528,2.7848) -- (4.4100,2.7600) -- (4.3876,2.7352) -- (4.3872,2.7120) -- (4.4088,2.6921) -- (4.4509,2.6769) -- (4.5106,2.6673) -- (4.5840,2.6640) -- (4.6659,2.6673) -- (4.7509,2.6769) -- (4.8330,2.6921) -- (4.9068,2.7120) -- (4.9672,2.7352) -- cycle;
\filldraw[fill=yellowY, draw=black!75, line width=0.3pt] (5.3846,2.5228) -- (5.3957,2.5012) -- (5.4282,2.4836) -- (5.4799,2.4711) -- (5.5473,2.4648) -- (5.6258,2.4649) -- (5.7100,2.4715) -- (5.7942,2.4841) -- (5.8727,2.5019) -- (5.9401,2.5237) -- (5.9918,2.5479) -- (6.0243,2.5730) -- (6.0354,2.5972) -- (6.0354,2.7972) -- (6.0243,2.8188) -- (5.9918,2.8364) -- (5.9401,2.8489) -- (5.8727,2.8552) -- (5.7942,2.8551) -- (5.7100,2.8485) -- (5.6258,2.8359) -- (5.5473,2.8181) -- (5.4799,2.7963) -- (5.4282,2.7721) -- (5.3957,2.7470) -- (5.3846,2.7228) -- cycle;
\filldraw[fill=yellowT, draw=black!75, line width=0.3pt] (6.0100,2.7600) -- (6.0324,2.7848) -- (6.0328,2.8080) -- (6.0112,2.8279) -- (5.9691,2.8431) -- (5.9094,2.8527) -- (5.8360,2.8560) -- (5.7541,2.8527) -- (5.6691,2.8431) -- (5.5870,2.8279) -- (5.5132,2.8080) -- (5.4528,2.7848) -- (5.4100,2.7600) -- (5.3876,2.7352) -- (5.3872,2.7120) -- (5.4088,2.6921) -- (5.4509,2.6769) -- (5.5106,2.6673) -- (5.5840,2.6640) -- (5.6659,2.6673) -- (5.7509,2.6769) -- (5.8330,2.6921) -- (5.9068,2.7120) -- (5.9672,2.7352) -- cycle;
\filldraw[fill=yellowY, draw=black!75, line width=0.3pt] (6.3846,2.5228) -- (6.3957,2.5012) -- (6.4282,2.4836) -- (6.4799,2.4711) -- (6.5473,2.4648) -- (6.6258,2.4649) -- (6.7100,2.4715) -- (6.7942,2.4841) -- (6.8727,2.5019) -- (6.9401,2.5237) -- (6.9918,2.5479) -- (7.0243,2.5730) -- (7.0354,2.5972) -- (7.0354,2.7972) -- (7.0243,2.8188) -- (6.9918,2.8364) -- (6.9401,2.8489) -- (6.8727,2.8552) -- (6.7942,2.8551) -- (6.7100,2.8485) -- (6.6258,2.8359) -- (6.5473,2.8181) -- (6.4799,2.7963) -- (6.4282,2.7721) -- (6.3957,2.7470) -- (6.3846,2.7228) -- cycle;
\filldraw[fill=yellowT, draw=black!75, line width=0.3pt] (7.0100,2.7600) -- (7.0324,2.7848) -- (7.0328,2.8080) -- (7.0112,2.8279) -- (6.9691,2.8431) -- (6.9094,2.8527) -- (6.8360,2.8560) -- (6.7541,2.8527) -- (6.6691,2.8431) -- (6.5870,2.8279) -- (6.5132,2.8080) -- (6.4528,2.7848) -- (6.4100,2.7600) -- (6.3876,2.7352) -- (6.3872,2.7120) -- (6.4088,2.6921) -- (6.4509,2.6769) -- (6.5106,2.6673) -- (6.5840,2.6640) -- (6.6659,2.6673) -- (6.7509,2.6769) -- (6.8330,2.6921) -- (6.9068,2.7120) -- (6.9672,2.7352) -- cycle;
\filldraw[fill=whiteX, draw=black!75, line width=0.5pt, line join=round] (9.1600,0.5600) -- (10.8400,1.8400) -- (10.8400,3.0400) -- (9.1600,1.7600) -- cycle;
\filldraw[fill=whiteY, draw=black!75, line width=0.5pt, line join=round] (7.1600,0.5600) -- (9.1600,0.5600) -- (9.1600,1.7600) -- (7.1600,1.7600) -- cycle;
\filldraw[fill=whiteT, draw=black!75, line width=0.5pt, line join=round] (7.1600,1.7600) -- (9.1600,1.7600) -- (10.8400,3.0400) -- (8.8400,3.0400) -- cycle;
\filldraw[fill=whiteY, draw=black!75, line width=0.3pt] (8.8046,2.8428) -- (8.8157,2.8212) -- (8.8482,2.8036) -- (8.8999,2.7911) -- (8.9673,2.7848) -- (9.0458,2.7849) -- (9.1300,2.7915) -- (9.2142,2.8041) -- (9.2927,2.8219) -- (9.3601,2.8437) -- (9.4118,2.8679) -- (9.4443,2.8930) -- (9.4554,2.9172) -- (9.4554,3.1172) -- (9.4443,3.1388) -- (9.4118,3.1564) -- (9.3601,3.1689) -- (9.2927,3.1752) -- (9.2142,3.1751) -- (9.1300,3.1685) -- (9.0458,3.1559) -- (8.9673,3.1381) -- (8.8999,3.1163) -- (8.8482,3.0921) -- (8.8157,3.0670) -- (8.8046,3.0428) -- cycle;
\filldraw[fill=whiteT, draw=black!75, line width=0.3pt] (9.4300,3.0800) -- (9.4524,3.1048) -- (9.4528,3.1280) -- (9.4312,3.1479) -- (9.3891,3.1631) -- (9.3294,3.1727) -- (9.2560,3.1760) -- (9.1741,3.1727) -- (9.0891,3.1631) -- (9.0070,3.1479) -- (8.9332,3.1280) -- (8.8728,3.1048) -- (8.8300,3.0800) -- (8.8076,3.0552) -- (8.8072,3.0320) -- (8.8288,3.0121) -- (8.8709,2.9969) -- (8.9306,2.9873) -- (9.0040,2.9840) -- (9.0859,2.9873) -- (9.1709,2.9969) -- (9.2530,3.0121) -- (9.3268,3.0320) -- (9.3872,3.0552) -- cycle;
\filldraw[fill=whiteY, draw=black!75, line width=0.3pt] (9.8046,2.8428) -- (9.8157,2.8212) -- (9.8482,2.8036) -- (9.8999,2.7911) -- (9.9673,2.7848) -- (10.0458,2.7849) -- (10.1300,2.7915) -- (10.2142,2.8041) -- (10.2927,2.8219) -- (10.3601,2.8437) -- (10.4118,2.8679) -- (10.4443,2.8930) -- (10.4554,2.9172) -- (10.4554,3.1172) -- (10.4443,3.1388) -- (10.4118,3.1564) -- (10.3601,3.1689) -- (10.2927,3.1752) -- (10.2142,3.1751) -- (10.1300,3.1685) -- (10.0458,3.1559) -- (9.9673,3.1381) -- (9.8999,3.1163) -- (9.8482,3.0921) -- (9.8157,3.0670) -- (9.8046,3.0428) -- cycle;
\filldraw[fill=whiteT, draw=black!75, line width=0.3pt] (10.4300,3.0800) -- (10.4524,3.1048) -- (10.4528,3.1280) -- (10.4312,3.1479) -- (10.3891,3.1631) -- (10.3294,3.1727) -- (10.2560,3.1760) -- (10.1741,3.1727) -- (10.0891,3.1631) -- (10.0070,3.1479) -- (9.9332,3.1280) -- (9.8728,3.1048) -- (9.8300,3.0800) -- (9.8076,3.0552) -- (9.8072,3.0320) -- (9.8288,3.0121) -- (9.8709,2.9969) -- (9.9306,2.9873) -- (10.0040,2.9840) -- (10.0859,2.9873) -- (10.1709,2.9969) -- (10.2530,3.0121) -- (10.3268,3.0320) -- (10.3872,3.0552) -- cycle;
\filldraw[fill=whiteY, draw=black!75, line width=0.3pt] (8.3846,2.5228) -- (8.3957,2.5012) -- (8.4282,2.4836) -- (8.4799,2.4711) -- (8.5473,2.4648) -- (8.6258,2.4649) -- (8.7100,2.4715) -- (8.7942,2.4841) -- (8.8727,2.5019) -- (8.9401,2.5237) -- (8.9918,2.5479) -- (9.0243,2.5730) -- (9.0354,2.5972) -- (9.0354,2.7972) -- (9.0243,2.8188) -- (8.9918,2.8364) -- (8.9401,2.8489) -- (8.8727,2.8552) -- (8.7942,2.8551) -- (8.7100,2.8485) -- (8.6258,2.8359) -- (8.5473,2.8181) -- (8.4799,2.7963) -- (8.4282,2.7721) -- (8.3957,2.7470) -- (8.3846,2.7228) -- cycle;
\filldraw[fill=whiteT, draw=black!75, line width=0.3pt] (9.0100,2.7600) -- (9.0324,2.7848) -- (9.0328,2.8080) -- (9.0112,2.8279) -- (8.9691,2.8431) -- (8.9094,2.8527) -- (8.8360,2.8560) -- (8.7541,2.8527) -- (8.6691,2.8431) -- (8.5870,2.8279) -- (8.5132,2.8080) -- (8.4528,2.7848) -- (8.4100,2.7600) -- (8.3876,2.7352) -- (8.3872,2.7120) -- (8.4088,2.6921) -- (8.4509,2.6769) -- (8.5106,2.6673) -- (8.5840,2.6640) -- (8.6659,2.6673) -- (8.7509,2.6769) -- (8.8330,2.6921) -- (8.9068,2.7120) -- (8.9672,2.7352) -- cycle;
\filldraw[fill=whiteY, draw=black!75, line width=0.3pt] (9.3846,2.5228) -- (9.3957,2.5012) -- (9.4282,2.4836) -- (9.4799,2.4711) -- (9.5473,2.4648) -- (9.6258,2.4649) -- (9.7100,2.4715) -- (9.7942,2.4841) -- (9.8727,2.5019) -- (9.9401,2.5237) -- (9.9918,2.5479) -- (10.0243,2.5730) -- (10.0354,2.5972) -- (10.0354,2.7972) -- (10.0243,2.8188) -- (9.9918,2.8364) -- (9.9401,2.8489) -- (9.8727,2.8552) -- (9.7942,2.8551) -- (9.7100,2.8485) -- (9.6258,2.8359) -- (9.5473,2.8181) -- (9.4799,2.7963) -- (9.4282,2.7721) -- (9.3957,2.7470) -- (9.3846,2.7228) -- cycle;
\filldraw[fill=whiteT, draw=black!75, line width=0.3pt] (10.0100,2.7600) -- (10.0324,2.7848) -- (10.0328,2.8080) -- (10.0112,2.8279) -- (9.9691,2.8431) -- (9.9094,2.8527) -- (9.8360,2.8560) -- (9.7541,2.8527) -- (9.6691,2.8431) -- (9.5870,2.8279) -- (9.5132,2.8080) -- (9.4528,2.7848) -- (9.4100,2.7600) -- (9.3876,2.7352) -- (9.3872,2.7120) -- (9.4088,2.6921) -- (9.4509,2.6769) -- (9.5106,2.6673) -- (9.5840,2.6640) -- (9.6659,2.6673) -- (9.7509,2.6769) -- (9.8330,2.6921) -- (9.9068,2.7120) -- (9.9672,2.7352) -- cycle;
\filldraw[fill=whiteY, draw=black!75, line width=0.3pt] (7.9646,2.2028) -- (7.9757,2.1812) -- (8.0082,2.1636) -- (8.0599,2.1511) -- (8.1273,2.1448) -- (8.2058,2.1449) -- (8.2900,2.1515) -- (8.3742,2.1641) -- (8.4527,2.1819) -- (8.5201,2.2037) -- (8.5718,2.2279) -- (8.6043,2.2530) -- (8.6154,2.2772) -- (8.6154,2.4772) -- (8.6043,2.4988) -- (8.5718,2.5164) -- (8.5201,2.5289) -- (8.4527,2.5352) -- (8.3742,2.5351) -- (8.2900,2.5285) -- (8.2058,2.5159) -- (8.1273,2.4981) -- (8.0599,2.4763) -- (8.0082,2.4521) -- (7.9757,2.4270) -- (7.9646,2.4028) -- cycle;
\filldraw[fill=whiteT, draw=black!75, line width=0.3pt] (8.5900,2.4400) -- (8.6124,2.4648) -- (8.6128,2.4880) -- (8.5912,2.5079) -- (8.5491,2.5231) -- (8.4894,2.5327) -- (8.4160,2.5360) -- (8.3341,2.5327) -- (8.2491,2.5231) -- (8.1670,2.5079) -- (8.0932,2.4880) -- (8.0328,2.4648) -- (7.9900,2.4400) -- (7.9676,2.4152) -- (7.9672,2.3920) -- (7.9888,2.3721) -- (8.0309,2.3569) -- (8.0906,2.3473) -- (8.1640,2.3440) -- (8.2459,2.3473) -- (8.3309,2.3569) -- (8.4130,2.3721) -- (8.4868,2.3920) -- (8.5472,2.4152) -- cycle;
\filldraw[fill=whiteY, draw=black!75, line width=0.3pt] (8.9646,2.2028) -- (8.9757,2.1812) -- (9.0082,2.1636) -- (9.0599,2.1511) -- (9.1273,2.1448) -- (9.2058,2.1449) -- (9.2900,2.1515) -- (9.3742,2.1641) -- (9.4527,2.1819) -- (9.5201,2.2037) -- (9.5718,2.2279) -- (9.6043,2.2530) -- (9.6154,2.2772) -- (9.6154,2.4772) -- (9.6043,2.4988) -- (9.5718,2.5164) -- (9.5201,2.5289) -- (9.4527,2.5352) -- (9.3742,2.5351) -- (9.2900,2.5285) -- (9.2058,2.5159) -- (9.1273,2.4981) -- (9.0599,2.4763) -- (9.0082,2.4521) -- (8.9757,2.4270) -- (8.9646,2.4028) -- cycle;
\filldraw[fill=whiteT, draw=black!75, line width=0.3pt] (9.5900,2.4400) -- (9.6124,2.4648) -- (9.6128,2.4880) -- (9.5912,2.5079) -- (9.5491,2.5231) -- (9.4894,2.5327) -- (9.4160,2.5360) -- (9.3341,2.5327) -- (9.2491,2.5231) -- (9.1670,2.5079) -- (9.0932,2.4880) -- (9.0328,2.4648) -- (8.9900,2.4400) -- (8.9676,2.4152) -- (8.9672,2.3920) -- (8.9888,2.3721) -- (9.0309,2.3569) -- (9.0906,2.3473) -- (9.1640,2.3440) -- (9.2459,2.3473) -- (9.3309,2.3569) -- (9.4130,2.3721) -- (9.4868,2.3920) -- (9.5472,2.4152) -- cycle;
\filldraw[fill=whiteY, draw=black!75, line width=0.3pt] (7.5446,1.8828) -- (7.5557,1.8612) -- (7.5882,1.8436) -- (7.6399,1.8311) -- (7.7073,1.8248) -- (7.7858,1.8249) -- (7.8700,1.8315) -- (7.9542,1.8441) -- (8.0327,1.8619) -- (8.1001,1.8837) -- (8.1518,1.9079) -- (8.1843,1.9330) -- (8.1954,1.9572) -- (8.1954,2.1572) -- (8.1843,2.1788) -- (8.1518,2.1964) -- (8.1001,2.2089) -- (8.0327,2.2152) -- (7.9542,2.2151) -- (7.8700,2.2085) -- (7.7858,2.1959) -- (7.7073,2.1781) -- (7.6399,2.1563) -- (7.5882,2.1321) -- (7.5557,2.1070) -- (7.5446,2.0828) -- cycle;
\filldraw[fill=whiteT, draw=black!75, line width=0.3pt] (8.1700,2.1200) -- (8.1924,2.1448) -- (8.1928,2.1680) -- (8.1712,2.1879) -- (8.1291,2.2031) -- (8.0694,2.2127) -- (7.9960,2.2160) -- (7.9141,2.2127) -- (7.8291,2.2031) -- (7.7470,2.1879) -- (7.6732,2.1680) -- (7.6128,2.1448) -- (7.5700,2.1200) -- (7.5476,2.0952) -- (7.5472,2.0720) -- (7.5688,2.0521) -- (7.6109,2.0369) -- (7.6706,2.0273) -- (7.7440,2.0240) -- (7.8259,2.0273) -- (7.9109,2.0369) -- (7.9930,2.0521) -- (8.0668,2.0720) -- (8.1272,2.0952) -- cycle;
\filldraw[fill=whiteY, draw=black!75, line width=0.3pt] (8.5446,1.8828) -- (8.5557,1.8612) -- (8.5882,1.8436) -- (8.6399,1.8311) -- (8.7073,1.8248) -- (8.7858,1.8249) -- (8.8700,1.8315) -- (8.9542,1.8441) -- (9.0327,1.8619) -- (9.1001,1.8837) -- (9.1518,1.9079) -- (9.1843,1.9330) -- (9.1954,1.9572) -- (9.1954,2.1572) -- (9.1843,2.1788) -- (9.1518,2.1964) -- (9.1001,2.2089) -- (9.0327,2.2152) -- (8.9542,2.2151) -- (8.8700,2.2085) -- (8.7858,2.1959) -- (8.7073,2.1781) -- (8.6399,2.1563) -- (8.5882,2.1321) -- (8.5557,2.1070) -- (8.5446,2.0828) -- cycle;
\filldraw[fill=whiteT, draw=black!75, line width=0.3pt] (9.1700,2.1200) -- (9.1924,2.1448) -- (9.1928,2.1680) -- (9.1712,2.1879) -- (9.1291,2.2031) -- (9.0694,2.2127) -- (8.9960,2.2160) -- (8.9141,2.2127) -- (8.8291,2.2031) -- (8.7470,2.1879) -- (8.6732,2.1680) -- (8.6128,2.1448) -- (8.5700,2.1200) -- (8.5476,2.0952) -- (8.5472,2.0720) -- (8.5688,2.0521) -- (8.6109,2.0369) -- (8.6706,2.0273) -- (8.7440,2.0240) -- (8.8259,2.0273) -- (8.9109,2.0369) -- (8.9930,2.0521) -- (9.0668,2.0720) -- (9.1272,2.0952) -- cycle;
\filldraw[fill=greenX, draw=black!75, line width=0.5pt, line join=round] (9.0000,2.4000) -- (9.8400,3.0400) -- (9.8400,4.2400) -- (9.0000,3.6000) -- cycle;
\filldraw[fill=greenY, draw=black!75, line width=0.5pt, line join=round] (5.0000,2.4000) -- (9.0000,2.4000) -- (9.0000,3.6000) -- (5.0000,3.6000) -- cycle;
\filldraw[fill=greenT, draw=black!75, line width=0.5pt, line join=round] (5.0000,3.6000) -- (9.0000,3.6000) -- (9.8400,4.2400) -- (5.8400,4.2400) -- cycle;
\filldraw[fill=greenY, draw=black!75, line width=0.3pt] (5.8046,4.0428) -- (5.8157,4.0212) -- (5.8482,4.0036) -- (5.8999,3.9911) -- (5.9673,3.9848) -- (6.0458,3.9849) -- (6.1300,3.9915) -- (6.2142,4.0041) -- (6.2927,4.0219) -- (6.3601,4.0437) -- (6.4118,4.0679) -- (6.4443,4.0930) -- (6.4554,4.1172) -- (6.4554,4.3172) -- (6.4443,4.3388) -- (6.4118,4.3564) -- (6.3601,4.3689) -- (6.2927,4.3752) -- (6.2142,4.3751) -- (6.1300,4.3685) -- (6.0458,4.3559) -- (5.9673,4.3381) -- (5.8999,4.3163) -- (5.8482,4.2921) -- (5.8157,4.2670) -- (5.8046,4.2428) -- cycle;
\filldraw[fill=greenT, draw=black!75, line width=0.3pt] (6.4300,4.2800) -- (6.4524,4.3048) -- (6.4528,4.3280) -- (6.4312,4.3479) -- (6.3891,4.3631) -- (6.3294,4.3727) -- (6.2560,4.3760) -- (6.1741,4.3727) -- (6.0891,4.3631) -- (6.0070,4.3479) -- (5.9332,4.3280) -- (5.8728,4.3048) -- (5.8300,4.2800) -- (5.8076,4.2552) -- (5.8072,4.2320) -- (5.8288,4.2121) -- (5.8709,4.1969) -- (5.9306,4.1873) -- (6.0040,4.1840) -- (6.0859,4.1873) -- (6.1709,4.1969) -- (6.2530,4.2121) -- (6.3268,4.2320) -- (6.3872,4.2552) -- cycle;
\filldraw[fill=greenY, draw=black!75, line width=0.3pt] (6.8046,4.0428) -- (6.8157,4.0212) -- (6.8482,4.0036) -- (6.8999,3.9911) -- (6.9673,3.9848) -- (7.0458,3.9849) -- (7.1300,3.9915) -- (7.2142,4.0041) -- (7.2927,4.0219) -- (7.3601,4.0437) -- (7.4118,4.0679) -- (7.4443,4.0930) -- (7.4554,4.1172) -- (7.4554,4.3172) -- (7.4443,4.3388) -- (7.4118,4.3564) -- (7.3601,4.3689) -- (7.2927,4.3752) -- (7.2142,4.3751) -- (7.1300,4.3685) -- (7.0458,4.3559) -- (6.9673,4.3381) -- (6.8999,4.3163) -- (6.8482,4.2921) -- (6.8157,4.2670) -- (6.8046,4.2428) -- cycle;
\filldraw[fill=greenT, draw=black!75, line width=0.3pt] (7.4300,4.2800) -- (7.4524,4.3048) -- (7.4528,4.3280) -- (7.4312,4.3479) -- (7.3891,4.3631) -- (7.3294,4.3727) -- (7.2560,4.3760) -- (7.1741,4.3727) -- (7.0891,4.3631) -- (7.0070,4.3479) -- (6.9332,4.3280) -- (6.8728,4.3048) -- (6.8300,4.2800) -- (6.8076,4.2552) -- (6.8072,4.2320) -- (6.8288,4.2121) -- (6.8709,4.1969) -- (6.9306,4.1873) -- (7.0040,4.1840) -- (7.0859,4.1873) -- (7.1709,4.1969) -- (7.2530,4.2121) -- (7.3268,4.2320) -- (7.3872,4.2552) -- cycle;
\filldraw[fill=greenY, draw=black!75, line width=0.3pt] (7.8046,4.0428) -- (7.8157,4.0212) -- (7.8482,4.0036) -- (7.8999,3.9911) -- (7.9673,3.9848) -- (8.0458,3.9849) -- (8.1300,3.9915) -- (8.2142,4.0041) -- (8.2927,4.0219) -- (8.3601,4.0437) -- (8.4118,4.0679) -- (8.4443,4.0930) -- (8.4554,4.1172) -- (8.4554,4.3172) -- (8.4443,4.3388) -- (8.4118,4.3564) -- (8.3601,4.3689) -- (8.2927,4.3752) -- (8.2142,4.3751) -- (8.1300,4.3685) -- (8.0458,4.3559) -- (7.9673,4.3381) -- (7.8999,4.3163) -- (7.8482,4.2921) -- (7.8157,4.2670) -- (7.8046,4.2428) -- cycle;
\filldraw[fill=greenT, draw=black!75, line width=0.3pt] (8.4300,4.2800) -- (8.4524,4.3048) -- (8.4528,4.3280) -- (8.4312,4.3479) -- (8.3891,4.3631) -- (8.3294,4.3727) -- (8.2560,4.3760) -- (8.1741,4.3727) -- (8.0891,4.3631) -- (8.0070,4.3479) -- (7.9332,4.3280) -- (7.8728,4.3048) -- (7.8300,4.2800) -- (7.8076,4.2552) -- (7.8072,4.2320) -- (7.8288,4.2121) -- (7.8709,4.1969) -- (7.9306,4.1873) -- (8.0040,4.1840) -- (8.0859,4.1873) -- (8.1709,4.1969) -- (8.2530,4.2121) -- (8.3268,4.2320) -- (8.3872,4.2552) -- cycle;
\filldraw[fill=greenY, draw=black!75, line width=0.3pt] (8.8046,4.0428) -- (8.8157,4.0212) -- (8.8482,4.0036) -- (8.8999,3.9911) -- (8.9673,3.9848) -- (9.0458,3.9849) -- (9.1300,3.9915) -- (9.2142,4.0041) -- (9.2927,4.0219) -- (9.3601,4.0437) -- (9.4118,4.0679) -- (9.4443,4.0930) -- (9.4554,4.1172) -- (9.4554,4.3172) -- (9.4443,4.3388) -- (9.4118,4.3564) -- (9.3601,4.3689) -- (9.2927,4.3752) -- (9.2142,4.3751) -- (9.1300,4.3685) -- (9.0458,4.3559) -- (8.9673,4.3381) -- (8.8999,4.3163) -- (8.8482,4.2921) -- (8.8157,4.2670) -- (8.8046,4.2428) -- cycle;
\filldraw[fill=greenT, draw=black!75, line width=0.3pt] (9.4300,4.2800) -- (9.4524,4.3048) -- (9.4528,4.3280) -- (9.4312,4.3479) -- (9.3891,4.3631) -- (9.3294,4.3727) -- (9.2560,4.3760) -- (9.1741,4.3727) -- (9.0891,4.3631) -- (9.0070,4.3479) -- (8.9332,4.3280) -- (8.8728,4.3048) -- (8.8300,4.2800) -- (8.8076,4.2552) -- (8.8072,4.2320) -- (8.8288,4.2121) -- (8.8709,4.1969) -- (8.9306,4.1873) -- (9.0040,4.1840) -- (9.0859,4.1873) -- (9.1709,4.1969) -- (9.2530,4.2121) -- (9.3268,4.2320) -- (9.3872,4.2552) -- cycle;
\filldraw[fill=greenY, draw=black!75, line width=0.3pt] (5.3846,3.7228) -- (5.3957,3.7012) -- (5.4282,3.6836) -- (5.4799,3.6711) -- (5.5473,3.6648) -- (5.6258,3.6649) -- (5.7100,3.6715) -- (5.7942,3.6841) -- (5.8727,3.7019) -- (5.9401,3.7237) -- (5.9918,3.7479) -- (6.0243,3.7730) -- (6.0354,3.7972) -- (6.0354,3.9972) -- (6.0243,4.0188) -- (5.9918,4.0364) -- (5.9401,4.0489) -- (5.8727,4.0552) -- (5.7942,4.0551) -- (5.7100,4.0485) -- (5.6258,4.0359) -- (5.5473,4.0181) -- (5.4799,3.9963) -- (5.4282,3.9721) -- (5.3957,3.9470) -- (5.3846,3.9228) -- cycle;
\filldraw[fill=greenT, draw=black!75, line width=0.3pt] (6.0100,3.9600) -- (6.0324,3.9848) -- (6.0328,4.0080) -- (6.0112,4.0279) -- (5.9691,4.0431) -- (5.9094,4.0527) -- (5.8360,4.0560) -- (5.7541,4.0527) -- (5.6691,4.0431) -- (5.5870,4.0279) -- (5.5132,4.0080) -- (5.4528,3.9848) -- (5.4100,3.9600) -- (5.3876,3.9352) -- (5.3872,3.9120) -- (5.4088,3.8921) -- (5.4509,3.8769) -- (5.5106,3.8673) -- (5.5840,3.8640) -- (5.6659,3.8673) -- (5.7509,3.8769) -- (5.8330,3.8921) -- (5.9068,3.9120) -- (5.9672,3.9352) -- cycle;
\filldraw[fill=greenY, draw=black!75, line width=0.3pt] (6.3846,3.7228) -- (6.3957,3.7012) -- (6.4282,3.6836) -- (6.4799,3.6711) -- (6.5473,3.6648) -- (6.6258,3.6649) -- (6.7100,3.6715) -- (6.7942,3.6841) -- (6.8727,3.7019) -- (6.9401,3.7237) -- (6.9918,3.7479) -- (7.0243,3.7730) -- (7.0354,3.7972) -- (7.0354,3.9972) -- (7.0243,4.0188) -- (6.9918,4.0364) -- (6.9401,4.0489) -- (6.8727,4.0552) -- (6.7942,4.0551) -- (6.7100,4.0485) -- (6.6258,4.0359) -- (6.5473,4.0181) -- (6.4799,3.9963) -- (6.4282,3.9721) -- (6.3957,3.9470) -- (6.3846,3.9228) -- cycle;
\filldraw[fill=greenT, draw=black!75, line width=0.3pt] (7.0100,3.9600) -- (7.0324,3.9848) -- (7.0328,4.0080) -- (7.0112,4.0279) -- (6.9691,4.0431) -- (6.9094,4.0527) -- (6.8360,4.0560) -- (6.7541,4.0527) -- (6.6691,4.0431) -- (6.5870,4.0279) -- (6.5132,4.0080) -- (6.4528,3.9848) -- (6.4100,3.9600) -- (6.3876,3.9352) -- (6.3872,3.9120) -- (6.4088,3.8921) -- (6.4509,3.8769) -- (6.5106,3.8673) -- (6.5840,3.8640) -- (6.6659,3.8673) -- (6.7509,3.8769) -- (6.8330,3.8921) -- (6.9068,3.9120) -- (6.9672,3.9352) -- cycle;
\filldraw[fill=greenY, draw=black!75, line width=0.3pt] (7.3846,3.7228) -- (7.3957,3.7012) -- (7.4282,3.6836) -- (7.4799,3.6711) -- (7.5473,3.6648) -- (7.6258,3.6649) -- (7.7100,3.6715) -- (7.7942,3.6841) -- (7.8727,3.7019) -- (7.9401,3.7237) -- (7.9918,3.7479) -- (8.0243,3.7730) -- (8.0354,3.7972) -- (8.0354,3.9972) -- (8.0243,4.0188) -- (7.9918,4.0364) -- (7.9401,4.0489) -- (7.8727,4.0552) -- (7.7942,4.0551) -- (7.7100,4.0485) -- (7.6258,4.0359) -- (7.5473,4.0181) -- (7.4799,3.9963) -- (7.4282,3.9721) -- (7.3957,3.9470) -- (7.3846,3.9228) -- cycle;
\filldraw[fill=greenT, draw=black!75, line width=0.3pt] (8.0100,3.9600) -- (8.0324,3.9848) -- (8.0328,4.0080) -- (8.0112,4.0279) -- (7.9691,4.0431) -- (7.9094,4.0527) -- (7.8360,4.0560) -- (7.7541,4.0527) -- (7.6691,4.0431) -- (7.5870,4.0279) -- (7.5132,4.0080) -- (7.4528,3.9848) -- (7.4100,3.9600) -- (7.3876,3.9352) -- (7.3872,3.9120) -- (7.4088,3.8921) -- (7.4509,3.8769) -- (7.5106,3.8673) -- (7.5840,3.8640) -- (7.6659,3.8673) -- (7.7509,3.8769) -- (7.8330,3.8921) -- (7.9068,3.9120) -- (7.9672,3.9352) -- cycle;
\filldraw[fill=greenY, draw=black!75, line width=0.3pt] (8.3846,3.7228) -- (8.3957,3.7012) -- (8.4282,3.6836) -- (8.4799,3.6711) -- (8.5473,3.6648) -- (8.6258,3.6649) -- (8.7100,3.6715) -- (8.7942,3.6841) -- (8.8727,3.7019) -- (8.9401,3.7237) -- (8.9918,3.7479) -- (9.0243,3.7730) -- (9.0354,3.7972) -- (9.0354,3.9972) -- (9.0243,4.0188) -- (8.9918,4.0364) -- (8.9401,4.0489) -- (8.8727,4.0552) -- (8.7942,4.0551) -- (8.7100,4.0485) -- (8.6258,4.0359) -- (8.5473,4.0181) -- (8.4799,3.9963) -- (8.4282,3.9721) -- (8.3957,3.9470) -- (8.3846,3.9228) -- cycle;
\filldraw[fill=greenT, draw=black!75, line width=0.3pt] (9.0100,3.9600) -- (9.0324,3.9848) -- (9.0328,4.0080) -- (9.0112,4.0279) -- (8.9691,4.0431) -- (8.9094,4.0527) -- (8.8360,4.0560) -- (8.7541,4.0527) -- (8.6691,4.0431) -- (8.5870,4.0279) -- (8.5132,4.0080) -- (8.4528,3.9848) -- (8.4100,3.9600) -- (8.3876,3.9352) -- (8.3872,3.9120) -- (8.4088,3.8921) -- (8.4509,3.8769) -- (8.5106,3.8673) -- (8.5840,3.8640) -- (8.6659,3.8673) -- (8.7509,3.8769) -- (8.8330,3.8921) -- (8.9068,3.9120) -- (8.9672,3.9352) -- cycle;
\node[font=\sffamily\bfseries\footnotesize, text=black, fill=white, fill opacity=0.9, text opacity=1, inner sep=1pt, rounded corners=1pt] at (2.0000,0.6000) {1};
\node[font=\sffamily\bfseries\footnotesize, text=black, fill=white, fill opacity=0.9, text opacity=1, inner sep=1pt, rounded corners=1pt] at (9.2000,0.6000) {2};
\node[font=\sffamily\bfseries\footnotesize, text=black, fill=white, fill opacity=0.9, text opacity=1, inner sep=1pt, rounded corners=1pt] at (5.0000,1.8000) {3};
\node[font=\sffamily\bfseries\footnotesize, text=black, fill=white, fill opacity=0.9, text opacity=1, inner sep=1pt, rounded corners=1pt] at (8.1600,1.1600) {4};
\node[font=\sffamily\bfseries\footnotesize, text=black, fill=white, fill opacity=0.9, text opacity=1, inner sep=1pt, rounded corners=1pt] at (7.0000,3.0000) {5};
\draw[line width=1.1pt, black!70] (13.0, 0.3) -- (14.5, 1.6);
\draw[line width=1.1pt, black!70] (16.0, 0.3) -- (14.5, 1.6);
\draw[line width=1.1pt, black!70] (16.0, 0.3) -- (17.5, 1.6);
\draw[line width=1.1pt, black!70] (14.5, 1.6) -- (16.0, 2.9);
\draw[line width=1.1pt, black!70] (17.5, 1.6) -- (16.0, 2.9);
\node[circle, draw=black!75, line width=0.6pt, fill=redT, minimum size=0.62cm, font=\sffamily\bfseries\small] at (13.0, 0.3) {1};
\node[circle, draw=black!75, line width=0.6pt, fill=blueT, minimum size=0.62cm, font=\sffamily\bfseries\small] at (16.0, 0.3) {2};
\node[circle, draw=black!75, line width=0.6pt, fill=yellowT, minimum size=0.62cm, font=\sffamily\bfseries\small] at (14.5, 1.6) {3};
\node[circle, draw=black!75, line width=0.6pt, fill=whiteT, minimum size=0.62cm, font=\sffamily\bfseries\small] at (17.5, 1.6) {4};
\node[circle, draw=black!75, line width=0.6pt, fill=greenT, minimum size=0.62cm, font=\sffamily\bfseries\small] at (16.0, 2.9) {5};
\end{tikzpicture}
\end{center}
\caption{A building counted in $a(5)$ and its contact graph.}\label{thebuilding}
\end{figure}
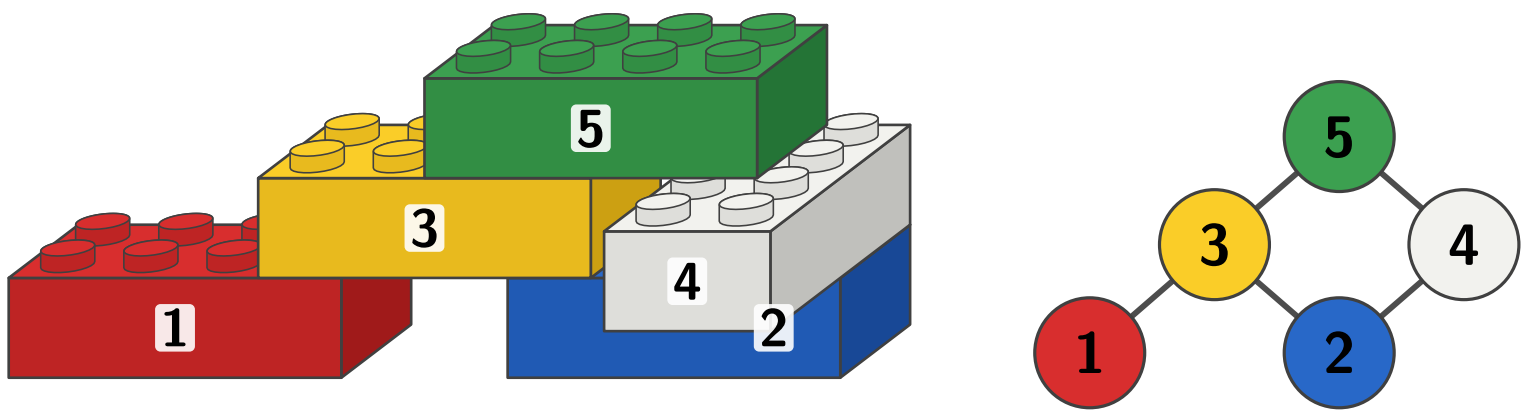
The known values of the sequence $a(n)$ are recorded in \Ano{112389}. Note that the count is not taking coloring of the bricks into account, and hence one may think of the buildings as being monochrome. We will, however, use colors in illustrations for visual clarity.

Obtaining $a(6)$ for the first time took the author several days on his laptop. Subsequent use of a mainframe computer 
and algorithmic improvements  with Abrahamsen \cite{asympenum} allowed 
computation of $a(7)$ and $a(8)$ in 2006, but predicted computing times measured in decades for $a(9)$.
In 2014, Johan Nilsson \cite{jn} determined this number
 using an embarrassingly parallel overcounting
method running part-time for about a month
on 200 cores.

The number $a(10)$ was obtained by
Matthias Simon \cite{ms} using an entirely
different approach allowing him to compute
the number on at most two standard laptops over
a period of 19 months. The
success of Simon's approach is
particularly striking in light of
the fact that it scales with the
number of connected simple graphs,
which grows superexponentially.

We have found
\[
a(11)=8364925094730088313
\]
using a development of the fundamental approach devised by Simon along with the use of HPC resources. The computation ran for less than 24 hours on GenomeDK at Aarhus University, through the DeiC national HPC collaboration, and required roughly 65,000 CPU hours.
 The new algorithm
also allows the recomputation of
$a(10)$ on a standard laptop (8 cores) in 39 wall clock hours, compared to Simon's
19 wall clock months.

The value of $a(11)$ was predicted by Gunning, Guttmann and Nilsson in \cite{agajgrmn:clc} to be 
\[
8.3649 \cdot
10^{18} \pm 2.4 \cdot10^{14}
\]
with the interval describing one standard deviation. Our exact computation confirms this, and is in fact only 0.1 standard deviations from the center of the predicted interval.

\subsection{Acknowledgments}

All of the computing for this project was performed on the GenomeDK cluster. We would like to thank GenomeDK and Aarhus University for providing computational resources and support that contributed to these research results.

We acknowledge substantial input from large language models, most prominently Anthropic's Claude on a Science Team Plan, both regarding the orchestration of the run and algorithmic improvements. 

\section{Simon's approach}

Simon's groundbreaking idea, as outlined in \cite{ms}, during a talk at the University of Copenhagen on November 13, 2018, is to
separate the dimensions in the
enumeration, taking into account
the differences between the $XY$-dimensions
on one hand and the $Z$-dimension on the other
when considering \lego bricks.
Working with labeled bricks,  and then correcting at the end,
Simon represents each \lego building
as a footprint in the $XY$-plane
combined with a layering describing
the $Z$ positions.

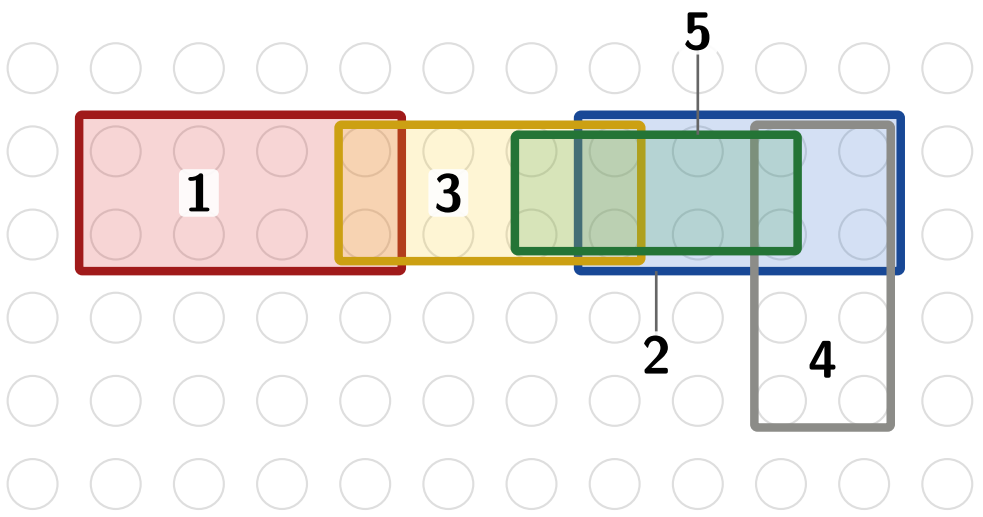
\begin{figure}
\begin{center}

\definecolor{blueT}{RGB}{40,104,200}
\definecolor{blueX}{RGB}{24,72,150}
\definecolor{blueY}{RGB}{32,90,180}
\definecolor{greenT}{RGB}{60,160,80}
\definecolor{greenX}{RGB}{36,116,54}
\definecolor{greenY}{RGB}{50,142,68}
\definecolor{outlineW}{RGB}{140,140,136}
\definecolor{redT}{RGB}{216,46,46}
\definecolor{redX}{RGB}{160,26,26}
\definecolor{redY}{RGB}{190,36,36}
\definecolor{whiteT}{RGB}{242,242,238}
\definecolor{whiteX}{RGB}{192,192,188}
\definecolor{whiteY}{RGB}{222,222,218}
\definecolor{yellowT}{RGB}{250,205,40}
\definecolor{yellowX}{RGB}{204,160,18}
\definecolor{yellowY}{RGB}{232,186,30}
\begin{tikzpicture}[scale=0.55]
\draw[black!13, line width=0.4pt] (-0.5,-2.5) circle (0.3);
\draw[black!13, line width=0.4pt] (-0.5,-1.5) circle (0.3);
\draw[black!13, line width=0.4pt] (-0.5,-0.5) circle (0.3);
\draw[black!13, line width=0.4pt] (-0.5,0.5) circle (0.3);
\draw[black!13, line width=0.4pt] (-0.5,1.5) circle (0.3);
\draw[black!13, line width=0.4pt] (-0.5,2.5) circle (0.3);
\draw[black!13, line width=0.4pt] (0.5,-2.5) circle (0.3);
\draw[black!13, line width=0.4pt] (0.5,-1.5) circle (0.3);
\draw[black!13, line width=0.4pt] (0.5,-0.5) circle (0.3);
\draw[black!13, line width=0.4pt] (0.5,0.5) circle (0.3);
\draw[black!13, line width=0.4pt] (0.5,1.5) circle (0.3);
\draw[black!13, line width=0.4pt] (0.5,2.5) circle (0.3);
\draw[black!13, line width=0.4pt] (1.5,-2.5) circle (0.3);
\draw[black!13, line width=0.4pt] (1.5,-1.5) circle (0.3);
\draw[black!13, line width=0.4pt] (1.5,-0.5) circle (0.3);
\draw[black!13, line width=0.4pt] (1.5,0.5) circle (0.3);
\draw[black!13, line width=0.4pt] (1.5,1.5) circle (0.3);
\draw[black!13, line width=0.4pt] (1.5,2.5) circle (0.3);
\draw[black!13, line width=0.4pt] (2.5,-2.5) circle (0.3);
\draw[black!13, line width=0.4pt] (2.5,-1.5) circle (0.3);
\draw[black!13, line width=0.4pt] (2.5,-0.5) circle (0.3);
\draw[black!13, line width=0.4pt] (2.5,0.5) circle (0.3);
\draw[black!13, line width=0.4pt] (2.5,1.5) circle (0.3);
\draw[black!13, line width=0.4pt] (2.5,2.5) circle (0.3);
\draw[black!13, line width=0.4pt] (3.5,-2.5) circle (0.3);
\draw[black!13, line width=0.4pt] (3.5,-1.5) circle (0.3);
\draw[black!13, line width=0.4pt] (3.5,-0.5) circle (0.3);
\draw[black!13, line width=0.4pt] (3.5,0.5) circle (0.3);
\draw[black!13, line width=0.4pt] (3.5,1.5) circle (0.3);
\draw[black!13, line width=0.4pt] (3.5,2.5) circle (0.3);
\draw[black!13, line width=0.4pt] (4.5,-2.5) circle (0.3);
\draw[black!13, line width=0.4pt] (4.5,-1.5) circle (0.3);
\draw[black!13, line width=0.4pt] (4.5,-0.5) circle (0.3);
\draw[black!13, line width=0.4pt] (4.5,0.5) circle (0.3);
\draw[black!13, line width=0.4pt] (4.5,1.5) circle (0.3);
\draw[black!13, line width=0.4pt] (4.5,2.5) circle (0.3);
\draw[black!13, line width=0.4pt] (5.5,-2.5) circle (0.3);
\draw[black!13, line width=0.4pt] (5.5,-1.5) circle (0.3);
\draw[black!13, line width=0.4pt] (5.5,-0.5) circle (0.3);
\draw[black!13, line width=0.4pt] (5.5,0.5) circle (0.3);
\draw[black!13, line width=0.4pt] (5.5,1.5) circle (0.3);
\draw[black!13, line width=0.4pt] (5.5,2.5) circle (0.3);
\draw[black!13, line width=0.4pt] (6.5,-2.5) circle (0.3);
\draw[black!13, line width=0.4pt] (6.5,-1.5) circle (0.3);
\draw[black!13, line width=0.4pt] (6.5,-0.5) circle (0.3);
\draw[black!13, line width=0.4pt] (6.5,0.5) circle (0.3);
\draw[black!13, line width=0.4pt] (6.5,1.5) circle (0.3);
\draw[black!13, line width=0.4pt] (6.5,2.5) circle (0.3);
\draw[black!13, line width=0.4pt] (7.5,-2.5) circle (0.3);
\draw[black!13, line width=0.4pt] (7.5,-1.5) circle (0.3);
\draw[black!13, line width=0.4pt] (7.5,-0.5) circle (0.3);
\draw[black!13, line width=0.4pt] (7.5,0.5) circle (0.3);
\draw[black!13, line width=0.4pt] (7.5,1.5) circle (0.3);
\draw[black!13, line width=0.4pt] (7.5,2.5) circle (0.3);
\draw[black!13, line width=0.4pt] (8.5,-2.5) circle (0.3);
\draw[black!13, line width=0.4pt] (8.5,-1.5) circle (0.3);
\draw[black!13, line width=0.4pt] (8.5,-0.5) circle (0.3);
\draw[black!13, line width=0.4pt] (8.5,0.5) circle (0.3);
\draw[black!13, line width=0.4pt] (8.5,1.5) circle (0.3);
\draw[black!13, line width=0.4pt] (8.5,2.5) circle (0.3);
\draw[black!13, line width=0.4pt] (9.5,-2.5) circle (0.3);
\draw[black!13, line width=0.4pt] (9.5,-1.5) circle (0.3);
\draw[black!13, line width=0.4pt] (9.5,-0.5) circle (0.3);
\draw[black!13, line width=0.4pt] (9.5,0.5) circle (0.3);
\draw[black!13, line width=0.4pt] (9.5,1.5) circle (0.3);
\draw[black!13, line width=0.4pt] (9.5,2.5) circle (0.3);
\draw[black!13, line width=0.4pt] (10.5,-2.5) circle (0.3);
\draw[black!13, line width=0.4pt] (10.5,-1.5) circle (0.3);
\draw[black!13, line width=0.4pt] (10.5,-0.5) circle (0.3);
\draw[black!13, line width=0.4pt] (10.5,0.5) circle (0.3);
\draw[black!13, line width=0.4pt] (10.5,1.5) circle (0.3);
\draw[black!13, line width=0.4pt] (10.5,2.5) circle (0.3);
\filldraw[fill=redT, fill opacity=0.2, draw=redX, line width=1.6pt, rounded corners=0.6pt] (0.06,0.06) rectangle (3.94,1.94);
\filldraw[fill=blueT, fill opacity=0.2, draw=blueX, line width=1.6pt, rounded corners=0.6pt] (6.06,0.06) rectangle (9.94,1.94);
\filldraw[fill=yellowT, fill opacity=0.2, draw=yellowX, line width=1.6pt, rounded corners=0.6pt] (3.18,0.18) rectangle (6.82,1.82);
\draw[draw=outlineW, line width=1.6pt, rounded corners=0.6pt] (8.18,-1.82) rectangle (9.82,1.82);
\filldraw[fill=greenT, fill opacity=0.2, draw=greenX, line width=1.6pt, rounded corners=0.6pt] (5.3,0.3) rectangle (8.7,1.7);
\node[font=\sffamily\bfseries\footnotesize, text=black, fill=white, fill opacity=0.9, text opacity=1, inner sep=1pt, rounded corners=1pt] at (1.5,1.0) {1};
\node[font=\sffamily\bfseries\footnotesize, text=black, fill=white, fill opacity=0.9, text opacity=1, inner sep=1pt, rounded corners=1pt] at (4.5,1.0) {3};
\node[font=\sffamily\bfseries\footnotesize, text=black, fill=white, fill opacity=0.9, text opacity=1, inner sep=1pt, rounded corners=1pt] at (9.0,-1.0) {4};
\draw[black!60, line width=0.5pt] (7.5,1.7) -- (7.5,2.75);
\node[font=\sffamily\bfseries\footnotesize, text=black, fill=white, fill opacity=0.9, text opacity=1, inner sep=1pt, rounded corners=1pt] at (7.5,2.95) {5};
\draw[black!60, line width=0.5pt] (7.0,0.06) -- (7.0,-0.75);
\node[font=\sffamily\bfseries\footnotesize, text=black, fill=white, fill opacity=0.9, text opacity=1, inner sep=1pt, rounded corners=1pt] at (7.0,-0.95) {2};
\end{tikzpicture}
\end{center}
\caption{The footprint of the building in Figure \ref{thebuilding}. The outlines of the bricks are shown with slightly decreasing sizes as their $Z$-coordinate increases.}\label{thefootprint}
\end{figure}

Inverting the count so as for it
to progress over graphs describing
how footprints overlap, Simon
uses inclusion-exclusion to eliminate
same-layer conditions of
non-overlapping and consider only
restrictions requiring certain bricks in the 2D footprints
to overlap. Reflecting the requirements
in Definition \ref{def:count}, these
graphs are all connected.

Finally, the observation that a
required overlap in the $XY$-plane
is described symmetrically in
$X$ and $Y$ allows the counting of
footprints to essentially take place
in one dimension.

We do not have access to
the concrete implementation
of Simon's algorithm, nor of the specifications of the hardware on which it ran,
but it can be inferred from \cite{ms} that the most expensive part  of the computation of $a(10)$ was the footprint count per graph. Indeed, Simon notes that the evaluation of each connected graph costs roughly 
 $4$ seconds per record, which multiplied by the number of connected graphs (cf.~\Ano{1349}) is  about $17.8$ months, a major part of the entire running time. Since the number of connected graphs with 11 vertices is about 86 times the number with 10,   a significant speedup of this aspect of Simon's method is necessary to make the computation of $a(11)$ practical, even employing HPC resources.
 
 \section{The footprint sweep}
In this section, we detail our algorithm for computing the number of labeled footprints (identified up to translation) realizing the overlap conditions from a given labeled connected graph $G$. Abbreviating $[n]=\{1,\dots,n\}$ to index both the vertices of $G$ and the bricks of the building, we will describe such a footprint by
\[
x,y:[n]\to \ZZ\qquad d:[n]\to \{H,V\}.
\]
With
\[
\ell^x(H)=4, \ell^x(V)=2, \ell^y(H)=2, \ell^y(V)=4
\]
the footprint of each brick $i\in[n]$ is
\[
[x(i),x(i)+\ell^x(d(i)))\times [y(i),y(i)+\ell^y(d(i))).
\]
Using Iverson notation, overlapping is indicated by
\[
o_{ij}=\Iver{\text{Bricks }i \text{ and } j\text{ overlap in positive area}}=o^x_{ij}o^y_{ij}
\]
where
\begin{eqnarray*}
o^x_{ij}&=&\Iver{x(i)<x(j)+\ell^x(d(j))} \Iver{x(j)<x(i)+\ell^x(d(i))}\\
o^y_{ij}&=&\Iver{y(i)<y(j)+\ell^y(d(j))} \Iver{y(j)<y(i)+\ell^y(d(i))}
\end{eqnarray*}

So more precisely:

\begin{probl}{pr:footprint}
Given a connected graph on $[n]$, compute the number $N(G)$ of $x,y,d$ with
\[
\min_{i\in[n]}x(i)=\min_{i\in[n]}y(i)=0
\]
and $o_{ij}=1$ whenever $ij$ is an edge of $G$.
\end{probl}

The symmetry of the overlap condition immediately allows the reduction 
\[
N(G)=\sum_{d:[n]\to\{H,V\}} X(G,d)X(G,\widetilde d)
\]
where $X(G,d)$ denotes the number of  possible choices of $x$ with $\min_{i\in[n]}x(i)=0$ so that $o^x_{ij}=1$ on all edges of $G$, and $\widetilde{d}$ reverses the symbols $H$ and $V$.
And 
when $d$ is given, a natural approach to compute the contribution to $N(G)$ is a \emph{deadline sweep} obtained by \emph{starting} the bricks in all $\emptyset\not= V_0\subseteq [n]$ by placing them at position $0$, and then imposing deadlines 
\[
\min_{i\in V_0,ij\in E(G)} \ell^x(d(i))-1
\]
on each  brick $j$ incident to one in $V_0$. Then at position $1$, all bricks incident to one placed vertically (i.e., with $d(i)=V$) must be placed, along with any number of bricks in $[n]\backslash V_0$, defining $V_1$. This in turn may impose additional deadlines, and by  connectivity of $G$, no brick will have  deadline later than $3n$, allowing  the number of placements in one direction to be computed as  number of different tuples
\[
(V_0,V_1,\dots,V_{3n}).
\]

The exponential character of this algorithm makes it impractical as stated, but as we shall explain shortly, it can be sped up significantly, especially at dense graphs, by merging states of an automaton describing the placement of the individual bricks. This will not suffice for our purposes, but we will explain this procedure first and present the one that allowed the computation to $n=11$ as a development of it. 

We describe the status of any brick by states $U$ (unplaced), $D_i$ (deadlines), $P$ (placed) and consider their life cycles as an automaton
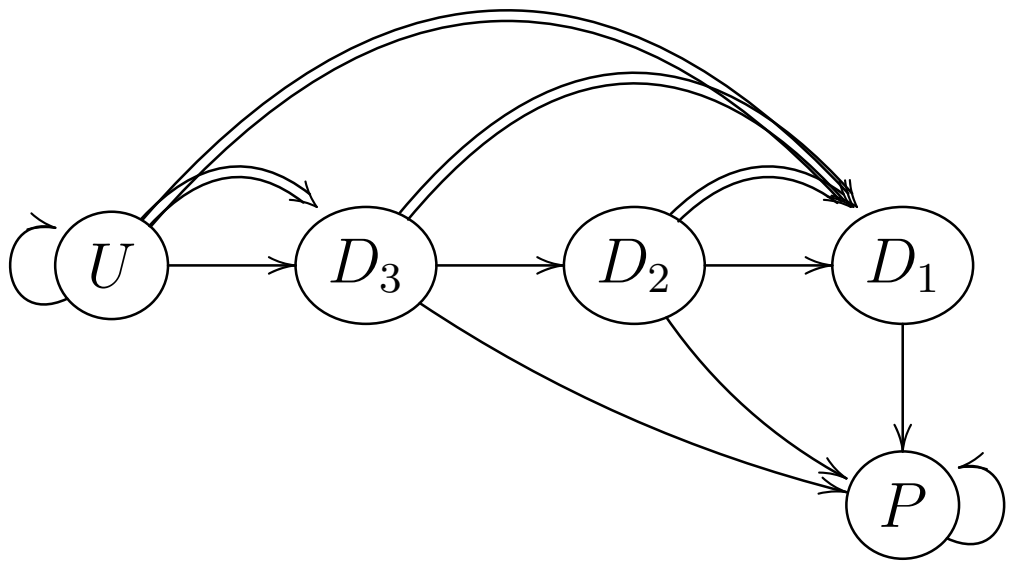
\[
\entrymodifiers={++[o][F-]}\xymatrix {
*{}&*{}&*{}&*{}&*{}\\
U\ar[r]\ar@(dl,ul)[]\ar@{=>}@/^1.5em/[r]\ar@{=>}@/^4em/[rrr]&D_3\ar[r]\ar@/_/[drr]\ar@{=>}@/^3em/[rr]&D_2\ar@{=>}@/^1.5em/[r]\ar[r]\ar@/_/[dr]&D_1\ar[d]\\
*{}&*{}&*{}&P\ar@(dr,ur)[]
}\]
The double arrows indicate steps that are forced by other bricks reaching state $P$ and then imposing deadlines on the one considered. When we ensure that at least one brick starts in state $P$, this will entail that the combined automaton ends up in the terminal state $(P,\dots,P)$ and that the minimal coordinate for bricks becomes $0$. Note that the placement of another brick can move the deadline to 1 even though a more generous deadline might already be given, depending on the direction vector $d$.

The key property of this description of the choice process is that it never records the coordinates of placed bricks, and consequently, states can be merged to reduce the computational space significantly. But because this must be done for all the $2^n$ choices of $d$, we were not able to bring the average time per graph at $n=11$ down below 4 seconds, which translates to at least 127 CPU years.

\newcommand{\Ps}{*+[o][F-]{P}}
\newcommand{\Us}{}
\newcommand{\Diiis}{*+[o][F-]{D_3}}
\newcommand{\Diis}{*+[o][F-]{D_2}}
\newcommand{\Dis}{*+[o][F-]{D_1}}
\newcommand{\Ds}{*+[o][F-]{D_0}}
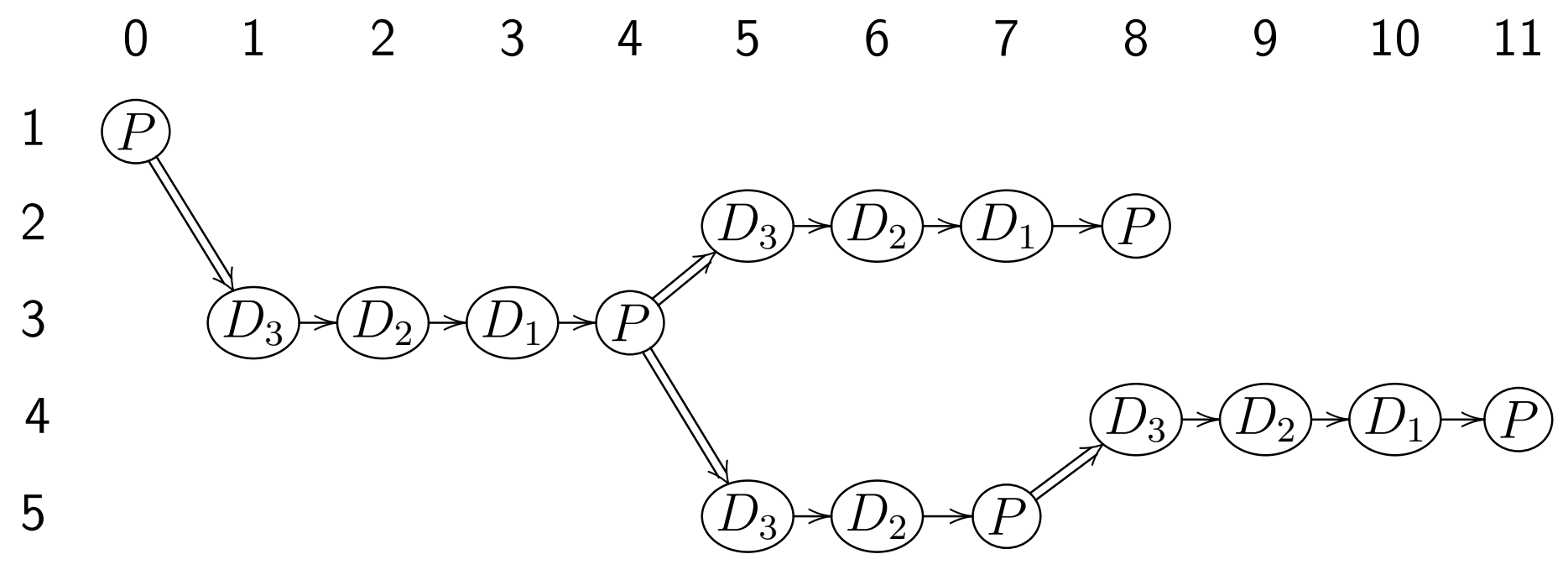
\begin{figure}
\[
\xymatrix@C=3mm@R=2mm{
&{\mathsf 0}&{\mathsf 1}&{\mathsf 2}&{\mathsf 3}&{\mathsf 4}&{\mathsf 5}&{\mathsf 6}&{\mathsf 7}&{\mathsf 8}&{\mathsf 9}&{\mathsf 1\mathsf 0}&{\mathsf 1\mathsf1}\\
{\mathsf 1}&\Ps\ar@{=>}[ddr]&\\
{\mathsf 2}&\Us&\Us&\Us&\Us&\Us&\Diiis\ar[r]&\Diis\ar[r]&\Dis\ar[r]&\Ps\\
{\mathsf 3}&\Us&\Diiis\ar[r]&\Diis\ar[r]&\Dis\ar[r]&\Ps\ar@{=>}[ur]\ar@{=>}[ddr]\\\
{\!\mathsf 4}&\Us&\Us&\Us&\Us&\Us&\Us&\Us&\Us&\Diiis\ar[r]&\Diis\ar[r]&\Dis\ar[r]&\Ps\\
{\mathsf 5}&\Us&\Us&\Us&\Us&\Us&\Diiis\ar[r]&\Diis\ar[r]&\Ps\ar@{=>}[ur]\\
}
\]
\caption{A sample run at $d=(H,H,H,H,V)$ leading to the $X$ coordinate of Figure \ref{thefootprint} (states $U$ and $P$ left out for legibility).}
\end{figure}

Our solution to this problem is to perform just one sweep for each $G$, allowing placements with deadline 3 throughout, but keeping track of which of the placed bricks would have to be placed in the $H$ direction to satisfy the overlap conditions encoded in $G$. This requires splitting the terminal state $P$ into two; $F$ for the situation where the brick is forced to be in the $H$ position, and $T$ (``turnable'') for the situation where it is not. We also introduce states $P_0$, $P_1$, $P_2$, $P_3$ to keep track of how long ago a brick was placed. If an edge $ij$ of $G$ is realized for $i$ already placed and in state $P_2$ or $P_3$, this forces the brick into the $F$ terminal state:

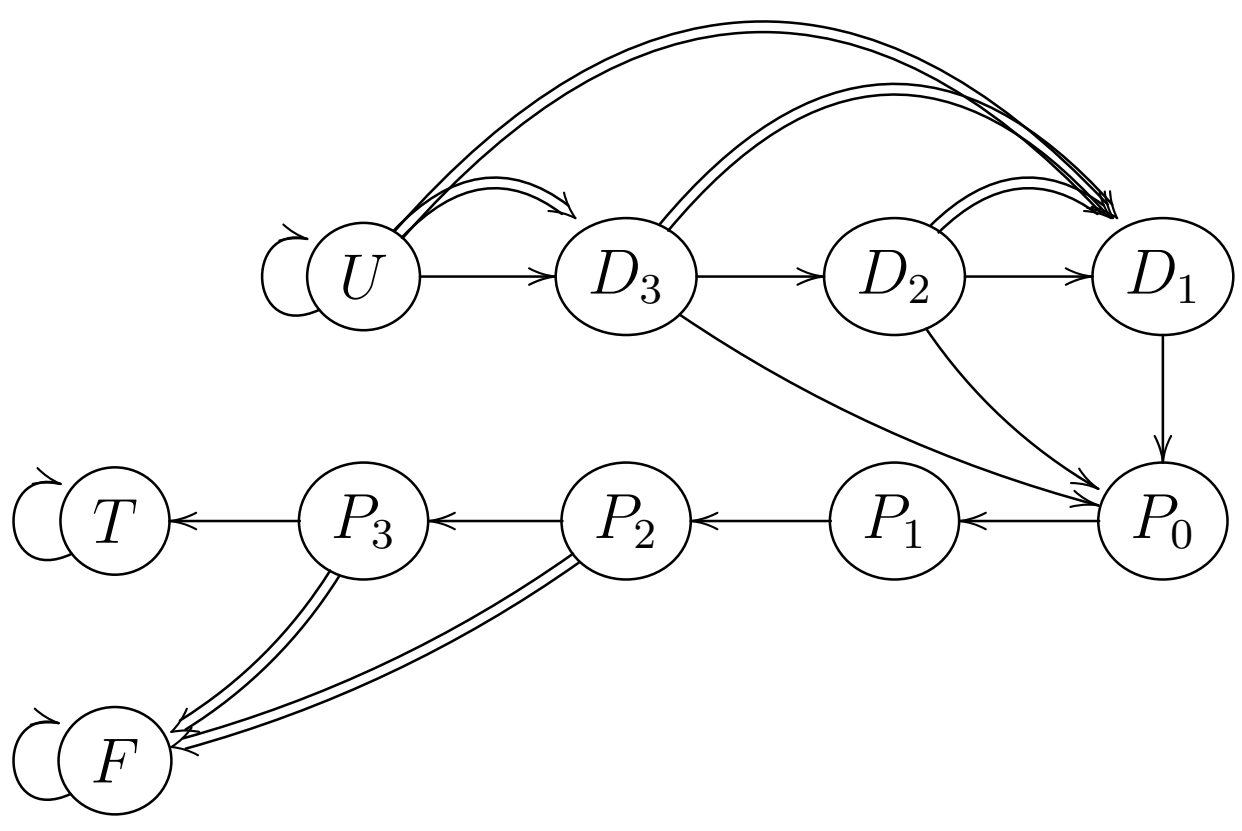
\[
\entrymodifiers={++[o][F-]}\xymatrix {
*{}&*{}&*{}&*{}&*{}\\
*{}&U\ar[r]\ar@(dl,ul)[]\ar@{=>}@/^1.5em/[r]\ar@{=>}@/^4em/[rrr]&D_3\ar[r]\ar@/_/[drr]\ar@{=>}@/^3em/[rr]&D_2\ar@{=>}@/^1.5em/[r]\ar[r]\ar@/_/[dr]&D_1\ar[d]\\
T\ar@(dl,ul)[]&P_3\ar[l]\ar@/^/@{=>}[dl]&P_2\ar[l]\ar@/^/@{=>}[dll]&P_1\ar[l]&P_0\ar[l]\\
F\ar@(dl,ul)[]&*{}&*{}&*{}&*{}
}\]

Using this setup for a single sweep allows the computation of $X(G,\mathcal F)$ counting the placements where exactly  $\mathcal F\subseteq [n]$ is forced into the $H$ direction. And then we use 

\begin{lemma}{cor:pairing}
\[
N(G)=
\sum_{\substack{{\mathcal L}_1, {\mathcal L}_2\subseteq [n]\\ {\mathcal L}_1\cap {\mathcal L}_2=\emptyset}} X(G,{\mathcal L}_1)X(G,{\mathcal L}_2)\,2^{\,n-|{\mathcal L}_1|-|{\mathcal L}_2|}.
\]
\end{lemma}
\begin{proof}
The first sweep forces bricks in $\mathcal L_1$ to be horizontal in the $x$ direction, while the second forces bricks in $\mathcal L_2$ to be horizontal in the $y$ direction. Thus the forced sets must be disjoint, since a brick cannot have both orientations. Each brick in neither set has two possible orientations.
\end{proof}

This approach brings the computation time per $G$ down to roughly $0.25$s on average, favoring dense graphs. The additional states of the new automaton increase computation mainly because less merging is possible, but this is more than amply made up for by not having to go over all choices of $d$. The overhead from the subset sum in Lemma \ref{cor:pairing} is negligible in this regard.

The state space is stored as four bit-planes over at most $11$ vertices, packed into a $44$-bit key. Subsets of freely chosen vertices  are enumerated with incrementally updated neighborhood unions. States at each position are merged by an LSD radix sort of the emitted (key, count)  pairs, saving significant time compared to  a hash table.

\section{Inverting the count}

We now detail how to turn the count of buildings into a sum over all connected graphs. Our approach is philosophically very close to Simon's but we organize the sum differently in order to speed up several aspects of the computation.

We denote by $F(n)$ the number of $n$-brick buildings up to translation in $\Z^3$. Since $a(n)$ counts them also up to rotation by multiples of $90^\circ$ about the vertical axis, we use Burnside's lemma, with $T_{180}(n)$ the number of half-turn-symmetric buildings counted up to planar rotation and $T_{90}(n)$ the number of quarter-turn-symmetric buildings counted up to planar rotation, to obtain
\[
a(n) \;=\; \frac{F(n) + 2T_{180}(n) + T_{90}(n)}{4} 
\]
It is clear that $T_{90}(n)=0$ unless $4\mid n$, and computing $T_{180}(n)$ can be done quickly for odd $n$ by noting that at least one brick must be rotated back to itself in any symmetric building, fixing the axis of rotation. We have found using two entirely different methods that 
\[
T_{180}(11)=2556043967
\]
(cf.~\Ano{123829}) and will not concern ourselves further with this aspect here.

We will again work with labeled buildings, and describe them as
\[
x,y,z:[n]\to \ZZ\qquad d:[n]\to \{H,V\}.
\]
 A non-zero translation never maps a finite set of bricks to itself, so the $n!$ labelings of a building are pairwise inequivalent and
\[
n!\,F(n) \;=\; \#\{\text{labeled buildings on } [n]\}\,/\,\text{translation}.
\]
It is convenient to fix a labeled building up to translation by imposing
\[
\min_{i\in[n]}x(i)=\min_{i\in [n]}y(i)=\min_{i\in [n]}z(i)=0
\]

To allow computations based on what it means to be a building, which is characterized as 
\begin{itemize}
\item no two bricks in the same layer overlap, and
\item the \emph{contact graph}, joining bricks in adjacent layers that overlap, is connected
\end{itemize}
we define 
\begin{eqnarray*}
A_z&=&\{ij\in \En\mid |z(i)-z(j)|=1\}\\
S_z&=&\{ij\in \En\mid z(i)=z(j)\}
\end{eqnarray*}
and
\[
G_{x,y,z,d}=\{ij\in A_z\mid o_{ij}=1\}
\]
where, as before, we have suppressed $x,y,d$ for the overlap indicators.
In Iverson notation, then
\begin{gather*}
\Iver{x,y,z,d\text{ defines a building}}=\\\prod_{ij\in S_z}\bigl(1-o_{ij}\bigr)\Iver{G_{x,y,z,d} \text{ is connected and spanning}}
\end{gather*}
where we can expand
\[
\Iver{G_{x,y,z,d} \text{ is connected and spanning}}= \sum_{\substack{H\subseteq A_z\\ H\ \text{conn.\ spanning}}}\ \prod_{ij\in H} o_{ij} \prod_{ij\in A_z\setminus H}\bigl(1-o_{ij}\bigr).
\]
and thus write the condition of being a building as
\begin{equation}\label{eq:building}
\prod_{ij\in S_z}\bigl(1-o_{ij}\bigr)\sum_{\substack{H\subseteq A_z\\ H\ \text{conn.\ spanning}}}\prod_{ij\in H} o_{ij} \prod_{ij\in A_z\setminus H}\bigl(1-o_{ij}\bigr)
\end{equation}
Because the $o_{ij}$ are idempotent, we can write this out
\[
\sum_{E\subseteq \En}c_z(E)\prod_{e\in E}o_e
\]
with the $c_z(E)$ to be determined.

\begin{lemma}{lem:whitney}
For $(x,y,z,d)$ as above, and $E\subseteq \En$, $c_z(E)$ does not depend on $x,y,d$ and is given by
\[
c_z(E)=\begin{cases}(-1)^{|E|}\kappa(E\cap A_z),& E\subseteq A_z\cup S_z,\\0,&\text{otherwise.}\end{cases}
\]
where
\[
\kappa(H) \;=\; \sum_{\substack{H'\subseteq H\\ H'\ \text{conn.\ spanning}}} (-1)^{|H'|}.
\]
\end{lemma}

Of course $\kappa(H)=0$ if $H$ is disconnected, and by work going back to Whitney \cite{greene1983interpretation,whitney}, $\kappa(H)$ can be related to more standard graph theoretical quantities in the connected case.
For instance  
\[
\kappa(H)=(-1)^{|V(H)|-1}\,T_H(1,0)\] where  $T_H$ is the Tutte polynomial, and $T_H(1,0)$ counts the acyclic orientations of $H$ with a single, prescribed source. Equivalently, $\kappa(H)$ is the M\"obius function of the bond lattice of $H$ between its bottom and top elements. This tells us that $\kappa$ is not cheaply computable, as the Tutte polynomial is known to be $\sharp P$-hard at $(1,0)$. We return to this problem in the next section.

\begin{proof}
It follows from \eqref{eq:building} that $c_z(E)=0$ unless $E\subseteq A_z\cup S_z$, and we abbreviate $E_A=E\cap A_z$, $E_S=E\cap S_z$. Expanding the terms of \eqref{eq:building} we have 
\[
\prod_{ij\in S_z}\bigl(1-o_{ij}\bigr) =\sum_{E_S\subseteq S_z}(-1)^{|E_S|}\prod_{e\in E_S}o_e
\]
and
\[
\sum_{H\subseteq E_A\subseteq A_z}(-1)^{|E_A\backslash H|}\prod_{e\in E_A}o_e
\]
for any fixed $H$ in the sum of \eqref{eq:building}, so the sign at every $E$ and $H$ becomes 
\[
(-1)^{|E\backslash H|}=(-1)^{|E|}(-1)^{|H|}.
\]
as desired.
\end{proof}

We note that $c_z(E)$ does only depend on the choice of $z$, as previewed by the notation.

Now we are ready to express the number of buildings as a sum over all labeled graphs $G$ on $[n]$:
\begin{eqnarray*}
n!F(n)&=&\sum_{x,y,z,d}\Iver{x,y,z,d\text{ defines a building}}/ \text{translation}\\
&=&\sum_{x,y,z,d}\sum_{G}c_z(G)\prod_{e\in E(G)}o_e/ \text{translation}\\
&=&\sum_{G}\sum_{\scriptsize \begin{array}{c}z\\\min z=0\end{array}}c_z(G)\sum_{\scriptsize \begin{array}{c}x,y,d\\\min x=0\\\min y=0\end{array}}\prod_{e\in E(G)}o_e\\
&=&\sum_{G}\sum_{z\in \Lambda(G)}(-1)^{|E(G)|}\kappa(G_z)\sum_{\scriptsize \begin{array}{c}x,y,d\\\min x=0\\\min y=0\end{array}}\prod_{e\in E(G)}o_e
\end{eqnarray*}
where
\[
E(G_z)=\{ij\in E(G)\mid |z(i)- z(j)|=1\}
\]
defines  the \emph{cut graph} of edges joining neighboring layers, and
\[
\Lambda(G)=\left\{z:[n]\to\ZZ\left|{\scriptsize\begin{array}{c}G_z\text{ connected}\\ij\in G\Longrightarrow |z(i)-z(j)|\leq 1\\ \min_{i\in[n]}z(i)=0\end{array}}\right.\right\}
\]
The definition of $\Lambda(G)$ reflects that by Lemma~\ref{lem:whitney}, a layering contributes only if its cut graph is connected, and that  terms in $c_z(E)$ contribute only when $E\subseteq A_z\cup S_z$.

Writing 
\[
C(G) \;=\; (-1)^{|E(G)|}\sum_{z\in\Lambda(G)} \kappa(G_z),
\]
we now have
\begin{theor}{thm:expansion}
Summing over labeled graphs $G$ on $[n]$,
\[
n!\,F(n) \;=\; \sum_{G} C(G)\,N(G),
\]
or in unlabeled form, \begin{equation}\label{eq:unlabeled}F(n)=\sum_{[G]} C(G)\,N(G)/|\Aut G|\end{equation} summing over isomorphism classes represented by a single $G$.
\end{theor}

As already discussed in the previous section, $N(G)$ counts placements modulo translation in which every edge of $G$ overlaps, and since $C(G)$ vanishes on non-connected $G$, the necessary connectivity can be assumed.

The sum is heavily signed. At $n=11$ it replaces about $3\times10^{19}$ buildings by about $10^{9}$ graphs, and the individual terms are many orders of magnitude larger than the result.

\section{Computing the coefficients \texorpdfstring{$C(G)$}{C(G)}}\label{sec:coefficients}

We must essentially compute $C(G)$  directly from Theorem~\ref{thm:expansion}, layering by layering, and when $z$ and thus $G_z$ is given, doing so is standard by a subset recursion over the vertices in $[n]$ that we shall outline below. 

We have, however, found a shortcut to compute the part of the sum stemming from  two-layer layerings that we shall describe here.
There are $2^n-2$ nonconstant two-level height assignments on $n$ vertices (though some are excluded from $\Lambda(G)$ because their cut graph is disconnected), and layerings with three or more layers become rare as the graphs get denser, so this provides a significant speedup in spite of the fact that there are relatively few two-level buildings.

It is straightforward to generate each layering in $\Lambda(G)$  exactly once using $z(1)=0$ as the translation normalization and then 
 giving each new vertex a height within $1$ of all of its already placed neighbors. This chooses the same layerings as the normalization $\min z=0$ in the definition of $\Lambda(G)$, although with different representatives. Removing two-layer layerings from the sum for separate analysis is also straightforward.
 
For any graph $H$ on $[n]$, we denote the subgraph induced at $S\subseteq [n]$ by $H[S]$: the vertices are $S$ and the edges those of $H$ that connect two vertices in $S$.
 We recall that
 \[
\kappa(H[S]) \;=\; \sum_{\substack{H'\subseteq H\\V(H')=S\\ H' \text{connected}}} (-1)^{|H'|}.
\]
 and relate this to
 \[
\iota(H[S])=\sum_{\substack{H'\subseteq H\\V(H')=S}} (-1)^{|H'|}
\]
in the following lemma.

\begin{lemma}{lem:subsetrecursion}
For any graph $H$ on $[n]$ and every nonempty $S\subseteq[n]$,
 \begin{equation}\label{eq:kapparec}
\kappa(H[S]) \;=\; \iota(H[S]) \;-\; \sum_{\substack{T\subsetneq S\\ \min S\,\in\, T}} \kappa(H[T])\,\iota(H[S\setminus T])
\end{equation}
 and
  \[
\iota(H[S]) \;=\; \Iver{H[S]\text{ has no edges}}\]
\end{lemma}
\begin{proof}
The description of $\iota$ follows directly from the binomial theorem. The recursive statement follows by grouping all terms in $\iota(H[S])$ according to the vertex set $T$ of the connected component which contains $\min S$. If $T=S$, the term goes into $\kappa(H[S])$, if not, the independence of the graphs on $T$ and $S\backslash T$ allows the multiplicative rewriting of the sum of the terms given.
\end{proof}

We can compute $\kappa(G_z)$ in 
 $3^{n-1}$ steps per layering, always having $1\in T$ and placing all other vertices either in $T$, in $S\backslash T$, or outside $S$.
 
 The  following result tempers the cost for many graphs.
 We write $\Lambda_2(G)$ for those layerings that have $z([n])=\{0,1\}$.

\begin{lemma}{twolayer}
For $n\ge2$ we have 
\[
\sum_{z\in\Lambda_2(G)} \kappa(G_z)=K^{[2]}([n]),
\]
where for all nonempty $S\subseteq [n]$ we set
\[
K^{[2]}(S)=\sum_{z:S\to \{0,1\}} \kappa(G_z[S]),
\]
and
\[
K^{[2]}(S) \;=\; 2^{\comp(G[S])} \;-\; \sum_{\substack{T\subsetneq S\\ \min S\,\in\, T}} K^{[2]}(T)\; 2^{\comp(G[S\setminus T])},
\]
where $\comp(\cdot)$ counts components. 
\end{lemma}
\begin{proof}
We can think of a two-layer layering as a $2$-coloring by $z$, where its cut graph is the set of bichromatic edges. Thus $K^{[2]}([n])$ agrees with the number we wish to compute because $\kappa$ vanishes on any non-connected $G_z$, in particular on the two constant colorings.

Since, in the recursion of \eqref{eq:kapparec},
the colorings of $T$ and $S\setminus T$ are independent, adding up all terms gives
\[
K^{[2]}(S) \;=\; I^{[2]}(S) \;-\; \sum_{\substack{T\subsetneq S\\ \min S\,\in\, T}} K^{[2]}(T)\,I^{[2]}(S\setminus T)
\]
where 
\[
I^{[2]}(S)=\sum_{z:S\to \{0,1\}} \iota(G_z[S]),
\]
Now, since
\begin{eqnarray*}
\iota(G_z[S\setminus T])&=&\Iver{G_z[S\setminus T]\text{ has no edges}}\\&=&\Iver{z\text{ is constant on all components of }G[S\setminus T]}
\end{eqnarray*}
the result follows since values either $0$ or $1$ can be assigned freely to each  component.
\end{proof}

\begin{remar}{rem:cutgraphs}
We realized after computing $a(11)$ that there are very few cut graphs contributing to $C(G)$. Indeed, they are bipartite and connected \Ano{5142}, and so there are only 25,598 relevant values of $\kappa$ needed at $n=11$. Even though there is substantial overhead from realizing isomorphisms and looking up precomputed values, this would have sped up our computation by at least 10\%.
\end{remar}

\section{Reduction to {2}-connected graphs}\label{sec:blocks}

As the only speedup mechanism recognizable from \cite{asympenum}, we now note that 
only $2$-connected graphs have to be evaluated. This is due to the fact that  the weight 
\[
w(G)=C(G)\,N(G)/2
\]
is multiplicative over blocks, allowing a standard  block decomposition reassembly  to reduce to this case. This not only gives a global reduction of the computation by the number of graphs following \Ano{2218} rather than \Ano{1349}, but restricts the computation to denser graphs, leading to speedups both concerning $N(G)$ and $C(G)$.
This is most valuable for $C(G)$, since for a typical $2$-connected graph on 11 vertices, two thirds of the contributing layerings have two layers.  The aggregation described in the previous section makes $C(G)$ $2.4$ times faster at $10$ vertices and $3.5$ times faster at 11 when restricted to $2$-connected graphs.

\begin{propo}{prop:mult}
Let $G$ be the union of $G_1$ and $G_2$, which share exactly one vertex $v$. Then $C(G)=C(G_1)\,C(G_2)$ and $N(G)/2=\bigl(N(G_1)/2\bigr)\bigl(N(G_2)/2\bigr)$.
\end{propo}

\begin{proof}
For $N$: once $v$'s orientation is fixed, the placements of the two sides relative to $v$ are independent, so $N(G)=\sum_{\omega_v}N_1(\omega_v)\,N_2(\omega_v)$. A quarter turn swaps the orientation of every brick and maps placements bijectively, so $N_i(H)=N_i(V)=N(G_i)/2$, and the product follows.

For $C$: a layering of $G$ is the same thing as a pair of layerings of $G_1$ and $G_2$ that agree at $v$, with heights taken relative to $v$. Its cut graph is the one-point join of the two cut graphs. A spanning subgraph of a one-point join is connected exactly when both halves are, so $\kappa$ is multiplicative, and so is the sign $(-1)^{|E|}$.
\end{proof}

\begin{theor}{thm:blocks} (Block composition.)
Let $b_k$ be the sum of $w(B)$ over labeled $2$-connected graphs $B$ on $[k]$, with $b_2=w(K_2)=92$ for the bridge, and $B(x)=\sum_k b_k x^k/k!$. The exponential generating function $R(x)$ of vertex-rooted connected graphs weighted by $w$ satisfies
\begin{equation}\label{eq:egf}
R(x) \;=\; x\,\exp\bigl(B'(R(x))\bigr),\qquad F(n) \;=\; \frac{2}{n}\,[x^n]\,R(x).
\end{equation}
\end{theor}

This is the standard decomposition of a rooted connected graph into the blocks at its root, each carrying rooted connected graphs at its other vertices, see e.g.~\cite{harary-palmer}, and holds for any weight that is multiplicative over blocks. Each connected graph has one more factor of $2$ in $C(G)N(G)$ than in the product of its block weights $w(B)=C(B)N(B)/2$, which explains the factor $2$ in the expression for $F(n)$.

We note that although $F(n)$ and all $b_k$ are integers, the recursion induced by \eqref{eq:egf} involves computations in $\QQ$ due to the term $2/n$ and the expansion of $\exp$. These are not time-critical and can be performed exactly. 

The block sums are graph-independent constants, and they differ wildly in sign and size:
\begin{center}
\begin{tabular}{|r|r|r|}
\hline
$k$ & $2$-connected graphs on $k$ vertices & $b_k$ \\
\hline\hline
2 & 1 & $92$ \\
3 & 1 & $-6936$ \\
4 & 3 & $933736$ \\
5 & 10 & $-149118480$ \\
6 & 56 & $7167748608$ \\
7 & 468 & $23395133971440$ \\
8 & 7123 & $-30743679807100800$ \\
9 & 194066 & $33428883328646077440$ \\
10 & 9743542 & $-36185834120433812640000$ \\
11 & 900969091 & 40651803567923756550693120\\\hline
\end{tabular}
\end{center}

\begin{remar}{rem:signs} Every graph's coefficient has a predictable sign. By Lemma~\ref{lem:whitney}, for a graph $G$ on $k$ vertices each $\kappa(G_z)$ is $(-1)^{k-1}$ times a count of acyclic orientations, so
\[
C(G) \;=\; (-1)^{\beta(G)}\,\lvert C(G)\rvert,\qquad \beta(G) \;=\; |E(G)| - k + 1,
\]
where $\beta(G)$ is the number of independent cycles of $G$. Since $N(G)>0$, $b_k$ is an alternating sum over $\beta$ of positive contributions that rise and fall smoothly, and these very nearly cancel:
\begin{center}
\begin{tabular}{|r|r|r|r|}\hline

$k$ & largest contribution & $b_k$ & cancellation \\
\hline\hline
5 & $-5.3\times10^{9}$ \ ($\beta=3$) & $-1.5\times10^{8}$ & $35\times$ \\
6 & $+1.6\times10^{13}$ \ ($\beta=4$) & $+7.2\times10^{9}$ & $2{,}200\times$ \\
7 & $+9.3\times10^{16}$ \ ($\beta=6$) & $+2.3\times10^{13}$ & $4{,}000\times$ \\
8 & $+1.1\times10^{21}$ \ ($\beta=8$) & $-3.1\times10^{16}$ & $36{,}000\times$ \\
\hline
\end{tabular}
\end{center}
The sign of $b_k$ is therefore a delicate remainder rather than a  mark of any one class of graphs: at $k=8$ the largest contribution is positive, yet the sum is negative. The observed pattern, $(-1)^k$ for $k\le6$ and $(-1)^{k+1}$ for $7\le k\le11$, has no explanation we know of.
\end{remar}

\section{Implementation and validation}\label{sec:implementation}

The computation is orchestrated  with a Python program \verb!run.py! delegating all heavy tasks either to a C++ program \verb!wcalc.cpp! or to  \verb!nauty! \cite{nauty} as follows:
\begin{enumerate}
\item \verb!run.py!: All 2-connected graphs are generated in the \verb!Graph6! format with \verb!nauty!'s \verb!geng! command (accessed via \verb!pynauty!), and then bundled into chunks that individual  threads will process in parallel. 
\item \verb!wcalc.cpp!: As many instances as possible  go through chunks of graphs and compute $C(G)$ and $N(G)$ for each, adding the numbers to a text file. This can be invoked as a batched job via Slurm.
\item \verb!run.py!: When all graphs of size $m$ have been processed, the $b_m$ are computed by summing \[
\frac{C(G)N(G)m!}{2|\Aut(G)|}
\]
The term $|\Aut(G)|$ is efficiently computed by \verb!nauty!. Then $F(n)$ is computed in exact rational arithmetic.
\end{enumerate}
\begin{figure}\begin{center}
\begin{verbatim}
                      DFw 2033792 70
                      DF{ 1789856 -72
                      DUW 2972752 -30
                      DUw 2122336 48
                      DU{ 1670048 -60
                      DV{ 1297504 84
                      D]w 1488760 -74
                      D]{ 1158368 90
                      D^{ 1000064 -112
                      D~{ 856352 150
\end{verbatim}\end{center}\caption{The full analysis at $n=5$. The first column contains the {\tt Graph6} description of the 10 2-connected graphs, and $N(G)$ and $C(G)$ are given in the next two.} 
\end{figure}

We ran the computation to $n=11$ on October 15--16, 2026,  using 112 nodes each of 32 CPUs. The collected CPU cost of the run, including pilots, testing, and archiving of data, was 65,362 hours.

The size of the computation does not allow an independent recomputation as things stand. We have in fact developed a rather fast GPU-based counter enumerating more directly as a development of the one in \cite{asympenum}. On the author's 2019 MacBook Pro, this program can recompute at $n=8$ in just under 200 seconds  (compared to our 440 hour computation in 2006), and at $n=9$ in 7.2 hours. But even taking into account the significant increase in GPU speed in recent years, we estimate a recomputation at $n=11$ to cost more than 5,000 GPU hours, and this  has not been available to us.

Thus we do not have a method-independent recomputation of $a(11)$. Our validation instead combines internal cross-checks with independent lower-order and $1\times2$ data. The main components are:

\begin{itemize}
\item \underline{Recomputation:} Using $F(1)=2$ and computing $b_2,\dots,b_{11}$, the program recomputes $a(1),\dots,a(10)$ after adjustment by $T_{180}(n)$ and $T_{90}(8)=122$. The numbers $a(1),\dots,a(9)$ have all been computed and recomputed in many different ways over the years, underpinning the correctness of $b_2,\dots,b_9$. And since these numbers in turn are signed sums of significantly larger summands, it is hard to imagine systematic or sporadic errors being cancelled out. The value $b_{10}$ gives a recalculation of Simon's value by completely  independent  software, but it must be admitted that the fundamental approaches are related. A different and independently coded recomputation of $a(10)$ by Lasse Deleuran in 2025 \cite{ld} also confirms the correctness of Simon's number.
\item \underline{Integrality:} As mentioned after Theorem \ref{thm:blocks}, the computation of $F(n)$ from $b_2,\dots,b_n$ is done in rationals by exact arithmetic, in a way that even though all the $b_k$ are integers, this is not automatically true for the resulting $F(n)$. The fact that all values computed are in fact integers supports their correctness, especially  at $n=11$ where the target integer was unknown.
\item \underline{Half size:}  With interval lengths $2$ and $1$ instead of $4$ and $2$, the same pipeline counts buildings of $1\times2$ bricks, for which the author produced in 2006 \Ano{123762} (all buildings $\widetilde{a}(n)$, up to rotation) and \Ano{123763} (the half-turn-symmetric ones $\widetilde{T}_{180}(n)$) using direct enumeration, which is sufficiently fast here. Whenever $4$ does not divide $n$, these give $\widetilde F(n)=4\widetilde a(n)-2\widetilde T_{180}(n)$, and the pipeline reproduces every such term through $n=11$. For  $n=8$ the implied number of quarter-turn-symmetric buildings comes out as $2$, which is easy to confirm. Because $C(G)$ does not depend on the brick dimensions, the $1\times2$ count at $11$ bricks reuses the $C(G)$ stored by the $2\times4$ computation, so it tests every one of those $900969091$ values, together with the graph list, the chunking and the composition, against independent data.
\item \underline{Alternative computations:} All four sweep counters give the same $N(G)$ as brute-force enumeration of placements on all $995$ connected graphs with $2$--$7$ vertices.
 The layering recursion, with and without two-layer aggregation, agrees with direct enumeration of layerings and brute-force counting of acyclic orientations on the same $995$ graphs to compute $C(G)$. The two versions also give identical output on all $7{,}123$ two-connected graphs on $8$ vertices.
\end{itemize}

\section{Outlook}

We end by briefly discussing the possibility of computing at $n=12$ and beyond. 

The author was convinced for a long time that a development of the basic enumeration program was the most promising option for reaching $n=11$, and hence put considerable effort (guided by LLMs) into optimizing such an approach for modern computer architecture. For this, he used his 2019 MacBook Pro with the ambition of passing to HPC for the last term.

But whereas the speedup was considerable, the fact that enumeration algorithms must to a large extent visit all buildings and hence are of a complexity reminiscent of the counts themselves was in the way of attempting an HPC run for $n=11$ based on such an approach. We also tried a transfer matrix method in light of the success of \cite{agajgrmn:clc}, but whereas this approach works very well in the 2D case, we could not see a viable way to our counting goal, and eventually turned to Simon's fundamental idea as already explained.

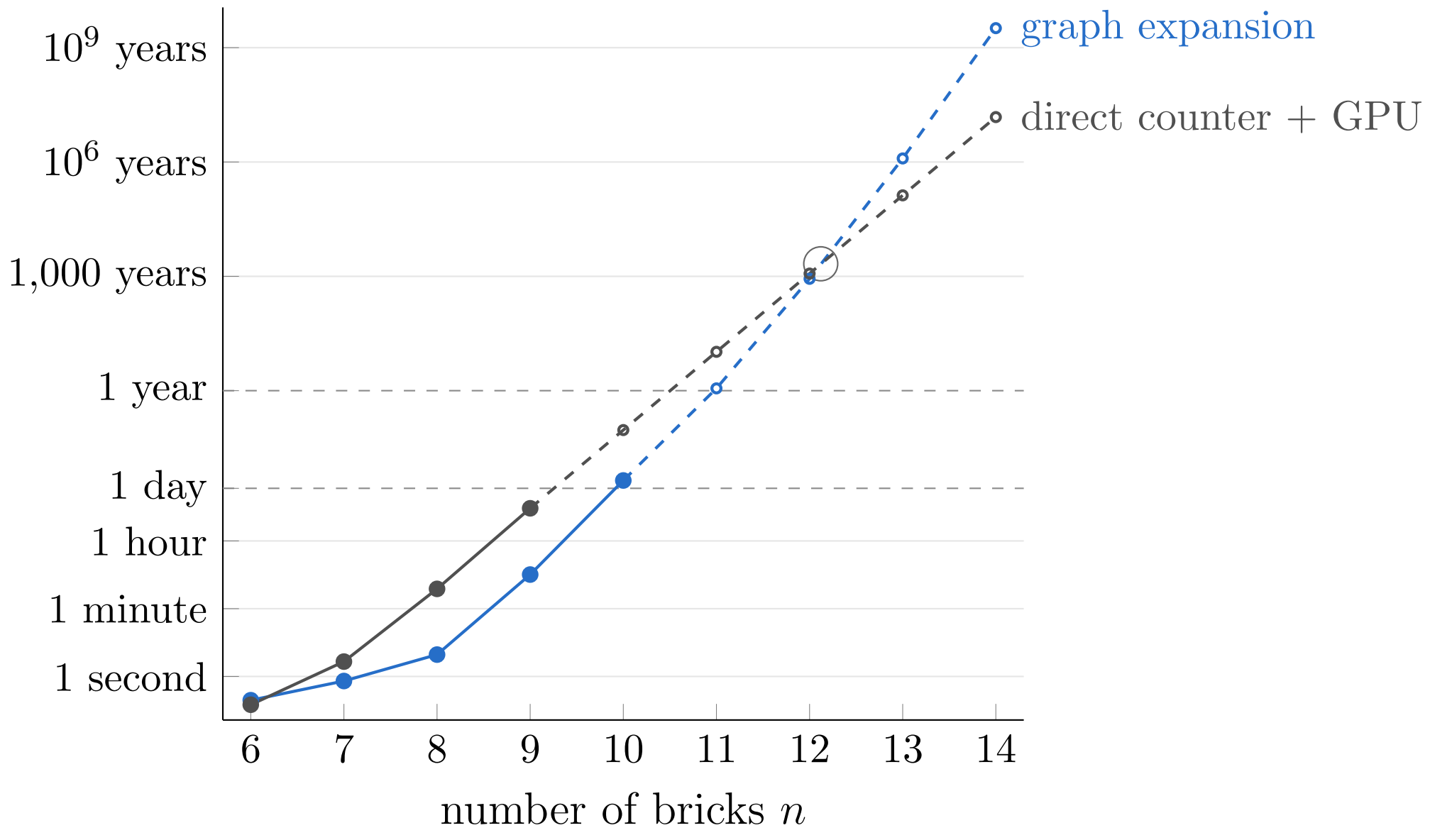
\begin{figure}
\centering
\begin{tikzpicture}
\begin{semilogyaxis}[
  width=0.66\textwidth, height=0.60\textwidth, clip=false,
  xmin=5.7, xmax=14.3, xtick={6,...,14}, xlabel={number of bricks $n$},
  ymin=2e-5, ymax=1e14,
  ytick={2.778e-4,1.667e-2,1,8.766e6,8.766e9,8.766e12},
  yticklabels={1 second,1 minute,1 hour,{1{,}000 years},$10^6$ years,$10^9$ years},
  extra y ticks={24,8766}, extra y tick labels={1 day,1 year},
  extra y tick style={grid=major, major grid style={dashed, black!45}},
  ymajorgrids, major grid style={black!10},
  tick label style={font=\small}, axis lines*=left]
\addplot[expcol, thick, mark=*, mark size=1.8pt] coordinates {(6,6.61e-5)(7,2.108e-4) (8,1.04e-3) (9,0.131) (10,38.54)};
\addplot[expcol, thick, dashed, mark=o, mark size=1.2pt, mark options={solid}] coordinates {(10,38.54)(11,1.007e4) (12,7.547e6) (13,1.087e10) (14,2.853e13)};
\addplot[dircol, thick, mark=*, mark size=1.8pt] coordinates {(6,5.031e-5) (7,6.817e-4) (8,0.05519) (9,7.178)};
\addplot[dircol, thick, dashed, mark=o, mark size=1.2pt, mark options={solid}] coordinates {(9,7.178) (10,811.1) (11,9.166e4) (12,1.036e7) (13,1.17e9) (14,1.323e11)};
\node[expcol, anchor=west, font=\small] at (axis cs:14.1,2.853e13) {graph expansion};
\node[dircol, anchor=west, font=\small] at (axis cs:14.1,1.323e11) {direct counter + GPU};
\node[circle, draw=black!60, minimum size=9pt, inner sep=0pt] (cross) at (axis cs:12.12,1.87e7) {};
\end{semilogyaxis}
\end{tikzpicture}
\caption{Wall time on the 2019 MacBook Pro. Filled marks are measured; hollow marks estimated. Under these extrapolations, the curves cross for the second time near $n=12$.}
\label{fig:mac}
\end{figure}

It remains obvious that the computing architecture and its future development significantly impacts how to aim for $n=12$ and beyond, but to fairly compare the enumeration approach to the algorithm presented here we have chosen to time our best implementations on the 2019 MacBook Pro as indicated in Figure \ref{fig:mac}, and extrapolate from there. Under these extrapolations, it shows first of all that $n=12$ and beyond is not an option with our code, but also that the window where Simon's approach is significantly faster than direct enumeration is about to close, because the number of connected graphs determining its complexity grows superexponentially.

We have a prototype for offloading a majority of the $N(G)$ computation to GPU, but whereas this would certainly lead to a speedup, such resources are presently out of reach with the current demand for GPU computing in the Danish and European HPC infrastructure.

It remains an unexplored possibility to divide the search space into parts more appropriate for certain algorithms and computing architectures, and work in parallel on them. A recent result by Deleuran \cite{ld} allowing the reduction of counts over layers with \emph{two} bricks rather than the single one used in \cite{asympenum} and Section \ref{sec:blocks} holds significant promise for such an approach.

\newcommand{\MR}[1]{}
\providecommand{\bysame}{\leavevmode\hbox to3em{\hrulefill}\thinspace}
\providecommand{\MR}{\relax\ifhmode\unskip\space\fi MR }
\providecommand{\MRhref}[2]{%
  \href{http://www.ams.org/mathscinet-getitem?mr=#1}{#2}
}
\providecommand{\href}[2]{#2}

\end{document}